\documentclass[11pt]{article}

\usepackage[margin=2.6cm]{geometry}
\usepackage{amsmath,amssymb,amsthm,mathtools}
\usepackage[round,authoryear]{natbib}
\usepackage[colorlinks,citecolor=blue,linkcolor=blue,urlcolor=blue]{hyperref}
\usepackage{graphicx,flafter,placeins,booktabs,array,enumitem,microtype}
\usepackage{xcolor}

\usepackage[utf8]{inputenc}
\usepackage[T1]{fontenc}
\usepackage[osf,sc]{mathpazo}

\newtheorem{theorem}{Theorem}[section]
\newtheorem{proposition}[theorem]{Proposition}

\newtheorem{corollary}[theorem]{Corollary}
\theoremstyle{definition}

\newcommand{\E}{\mathbb E}

\title{The Random Analyst}
\author{%
  G\'erard Biau\\[4pt]
  \small Sorbonne Universit\'e, Institut universitaire de France\\[2pt]
  \small \texttt{gerard.biau@sorbonne-universite.fr}
}
\date{}

\begin{document}
\maketitle

\begin{abstract}
\noindent
Modern data-analysis pipelines increasingly involve searching among candidate analyses before reporting an estimate. We investigate the replicability of the resulting conclusions by studying independent repetitions of the entire search-and-report procedure on new samples. Can the reported estimates become increasingly similar across studies while remaining systematically displaced from the scientific target? We first show that this can happen in a simple Gaussian model. Although every prespecified candidate estimator is unbiased and consistent, searching over exponentially many candidates within each study yields a selected estimator converging to a limit different from the target, a selection-induced pseudotruth. Across independent repetitions, the reported
estimates thus become replicable but systematically
inaccurate, while the selected candidates differ
with probability tending to one. We then study dependent Gaussian candidates, including marginal regression coefficients, and show that Gaussian width governs the displacement, which persists when squared width is of the same order as sample size. Metric entropy and generic chaining reveal that the effective complexity of the search is determined by the geometry of statistically distinguishable analyses rather than by their nominal number. Large-deviation arguments further extend the
phenomenon beyond Gaussian models. Numerical agreement between independent
studies therefore need not indicate either a stable
analytical path or proximity to the scientific target.
\smallskip

\end{abstract}

\section{Introduction}
\label{sec:introduction}

Statistical procedures routinely use the data to decide how the data should be
analyzed. Variables may be selected, models compared, tuning parameters chosen,
transformations tried, and hypotheses screened before a final result is
reported. There is nothing unusual about such data dependence: cross-validation,
model selection, and adaptive estimation are themselves standard procedures,
and the complete mapping from the data to the reported conclusion may be fully
specified before any observations are seen.

It is also classical that selection can distort subsequent inference. When the same data are used both to search among candidate analyses and to report the selected result, favorable sampling fluctuations that influence selection may be carried into the reported conclusion. This mechanism underlies the winner's curse
\citep{andrews2024} and more generally complicates post-selection
inference \citep{berk2013,russo2016,mclatchie2024}. Methods for valid inference in this setting
include sample splitting, simultaneous inference, and
conditional selective inference; see
\citet{fithian2014} and the review of
\citet{kuchibhotla2022}.

We study a different repetition experiment. Rather than conditioning on one
realized selection, we independently repeat the \emph{entire search-and-report
procedure} on new data. Each study generates a new sample, searches again over
the available analyses, possibly selects a different one, and reports the
corresponding conclusion. Our question is whether those reported conclusions
can become increasingly similar across studies while remaining systematically
displaced from the scientific target.

A simple decomposition makes the question precise. Let \(D_n\) denote the data, \(\mathcal A_n\) a family of candidate analyses, and
\[
A_n=A_n(D_n)\in\mathcal A_n
\]
the analysis selected from that family using \(D_n\) (for simplicity, we consider a deterministic selection rule here). If \(T_n(D_n,a)\) is the conclusion associated with candidate \(a\), the reported conclusion is
\[
T_n^\star=T_n\bigl(D_n,A_n(D_n)\bigr).
\]
Thus the same data enter the procedure twice: they determine which analysis is selected and they are also used to compute the conclusion reported by that selected analysis. Even when the selection rule
$D_n \mapsto A_n(D_n)$ is fully prespecified,
the identity of the selected analysis
remains data-dependent. Our benchmark is that every deterministic sequence \(a_n\in\mathcal A_n\), fixed independently of the realized data, satisfies
\[
T_n(D_n,a_n)\xrightarrow{P}\theta.
\]
Thus any persistent displacement of \(T_n^\star\) is generated by the data-dependent search over the family rather than by the asymptotic behavior of a fixed candidate. This is a modeling benchmark, not a claim that every practical collection of analyses has a common target. One setting in which the common-target benchmark is natural is covariate adjustment in randomized experiments: under suitable conditions, different regression adjustments can consistently estimate the same average treatment effect \citep{lin2013}.

For two independent repetitions, write \(T_n^{\star(1)}\) and
\(T_n^{\star(2)}\) for the reported conclusions and \(A_n^{(1)}\),
\(A_n^{(2)}\) for the selected analytical paths. We distinguish target
consistency, conclusion replicability, and path replicability. The central phenomenon is that
\[
T_n^{\star(1)}-T_n^{\star(2)}\xrightarrow{P}0
\]can hold even though \(T_n^\star\not\xrightarrow{P}\theta\). In particular, the reported conclusion may instead converge to a displaced limit,
\[
T_n^\star\xrightarrow{P}\theta+b,\qquad b\neq0,
\]even while the selected paths become increasingly unlikely to coincide.

The simplest example selects the largest among \(L_n\) independent unbiased
Gaussian estimators with standard error \(\sigma/\sqrt n\). If $\log L_n/n\to c>0$, the selected estimator
converges not to $\theta$ but to
$\theta+\sigma\sqrt{2c}$. Even without this growth condition,
the reported conclusions in two independent studies become
asymptotically indistinguishable, whereas the maximizing indices
agree with probability only $1/L_n$. This example is only a
benchmark. The main result of the paper is the correlated Gaussian analogue.
If \(w_n\) denotes the Gaussian width of the candidate family, then
\[
T_n^\star
=
\theta+\frac{\sigma w_n}{\sqrt n}+O_P(n^{-1/2}).
\]
Here \(w_n\) measures the average displacement generated by search, in units of the reference standard error, while fluctuations around this displaced center remain at most on the ordinary scale. The critical scale is therefore
\(w_n^2\) of order \(n\). Gaussian width replaces nominal candidate count by
the geometry of the joint sampling fluctuations; the canonical Gaussian
metric, packing and covering numbers, and generic chaining make this geometry
intrinsic.

The distinction from ordinary holdout evaluation is useful. If data \(D_n\)
select \(A_n\) but the selected analysis is evaluated on genuinely independent
data \(D_n'\), then the fluctuations driving selection do not enter the final
evaluation. Under the common-target benchmark, target consistency is restored.
The issue studied in the present paper is therefore not data dependence alone, but selection aligned with sampling fluctuations that also affect the reported conclusion.

Our perspective is related to several established literatures. Adaptive data analysis studies how repeated data-dependent querying affects statistical validity and generalization, as well as the bias induced by adaptive selection \citep{dwork2015a,dwork2015b,russo2016}. Post-selection inference studies the validity and limitations of inference after a data-dependent choice \citep{leebpotscher2005,leebpotscher2006,leebpotscher2008}. Proposed approaches include conditional selective inference \citep{fithian2014,taylor2015} and simultaneous protection over possible selected models, as in PoSI \citep{berk2013,bachoc2020}. These approaches address the validity of inference following selection. Our object is instead the unconditional law obtained when the complete search-and-report procedure is rerun on independent data. Reproducible-learning algorithms make stability across independent samples a design objective \citep{impagliazzo2022}. By contrast, we study settings in which stability can arise in the reported conclusion without stability of the selected path.

Researcher degrees of freedom and specification search long predate automated
analysis \citep{white2000,simmons2011,gelmanloken2013}, and many-analyst studies document
substantial analytical variation on a common dataset
\citep{silberzahn2018,botvinik2020,breznau2022}. Modern automated and agentic
systems can greatly increase the number and diversity of pipelines explored at
low cost; \citet{bertran2026} illustrates this growing search capacity on a
fixed dataset. Although automation is not required for the statistical
mechanism studied here, it makes the large-search regimes in which the
mechanism becomes pronounced increasingly plausible in practice.

\paragraph{Organization of the paper.}
Section~\ref{sec:framework} introduces the replication law and distinguishes target consistency, conclusion replicability, and path replicability. We then turn in Section~\ref{sec:minimal} to the independent Gaussian benchmark, which isolates the basic mechanism and its implications for superconcentration, directional selection, and proxy-based search. Section~\ref{sec:correlated} moves beyond independence and establishes the main Gaussian-width theorem. This geometric viewpoint is developed further in Section~\ref{sec:search-dimension} through correlated examples and intrinsic Gaussian-process geometry. Section~\ref{sec:large-deviations} shows how the same mechanism extends beyond Gaussianity, while Section~\ref{sec:scope} discusses diagnostics, correction, and limitations. Proofs and technical derivations omitted from the main text are collected in the Supplementary Material.

\section{Replication of a search-and-report procedure}
\label{sec:framework}

Let \(D_n=(X_1,\ldots,X_n)\sim P^n\) be a sample, where the observations take
values in a measurable space \(\mathcal X\), and let \(\theta=\theta(P)\) be a
scalar scientific target. Let \(\mathcal A_n\) be a measurable space of
candidate analyses and
\[
T_n:\mathcal X^n\times\mathcal A_n\to\mathbb R
\]
a measurable map, with \(T_n(D_n,a)\) denoting the conclusion associated with
candidate \(a\in\mathcal A_n\).

The selected analysis is a random element
$A_n\in\mathcal A_n$.
We refer to the particular candidate ultimately
selected by the procedure as the \emph{analytical path}.
The randomness of this selection is
the sense in which the analyst of the title
is random. To make the dependence on the data explicit,
a deterministic selection rule can be written
\[
A_n=A_n(D_n)=\Phi_n(D_n).
\]
More generally, a randomized selection rule is described by
\[
A_n\mid D_n\sim Q_n(\cdot\mid D_n),
\]
where \(Q_n\) is a probability kernel. The kernel allows internal randomization,
although none of the main mechanisms below requires it. In either case, the
reported conclusion is
\[
T_n^\star=T_n(D_n,A_n).
\]

\paragraph{The replication law.}
Holding a realized dataset \(D_n=x\) fixed and rerunning a stochastic selector
reveals \(\mathcal L(T_n^\star\mid D_n=x)\). Many-analyst experiments on a
common dataset, including \citet{bertran2026}, describe variation of this
conditional kind. They reveal how conclusions vary across analyses of one
realized dataset, but do not by themselves describe what happens when the
sampling experiment itself is repeated. Our object instead regenerates both
the data and the analysis selected in response to them. Averaging over repeated
datasets gives
\[
\mathcal L(T_n^\star)
=
\int \mathcal L(T_n^\star\mid D_n=x)\,P^n(\mathrm dx),
\]
which we call the \emph{replication law}. When \(A_n=a_n\) is fixed before observing the
data, this reduces to the ordinary sampling law of \(T_n(D_n,a_n)\).
Thus the replication law integrates both sources of variation: fresh sampling and
the analysis selected in response to the resulting data.

For two independent repetitions, let \((D_n^{(1)},A_n^{(1)})\) and
\((D_n^{(2)},A_n^{(2)})\) be independent copies generated by the same complete
procedure and write
\[
T_n^{\star(r)}=T_n(D_n^{(r)},A_n^{(r)}),\qquad r=1,2.
\]
Each repetition draws new data and reruns the search. This is replicability
across independently obtained data, rather than computational reproducibility
on a fixed dataset \citep{nasem2019}.

\paragraph{Target, conclusion, and path.}
We call \(T_n^\star\xrightarrow{P}\theta\) \emph{target consistency}; it asks
whether the complete search-and-report procedure approaches the intended
scientific target. We call
\(
T_n^{\star(1)}-T_n^{\star(2)}\xrightarrow{P}0
\)
\emph{conclusion replicability}; it requires the reported conclusions from independent repetitions to become asymptotically indistinguishable, without requiring convergence to any fixed value. Finally, for a discrete candidate family, path agreement
can be measured by \(\mathbb P(A_n^{(1)}=A_n^{(2)})\); for a continuous family,
exact equality is generally too strong and a suitable distance between the
selected analyses is more informative. Conclusion and path replicability are
therefore distinct properties of the same repeated experiment.

\paragraph{The common-target benchmark.}
Throughout the paper we impose the common-target condition that, for every
deterministic sequence \(a_n\in\mathcal A_n\) fixed independently of the
realized data,
\[
T_n(D_n,a_n)\xrightarrow{P}\theta.
\]
This is deliberately stronger than requiring consistency of one particular
candidate. It rules out, in particular, a persistently biased analysis hidden
inside \(\mathcal A_n\): no deterministic choice from the candidate family can
generate an asymptotic displacement from \(\theta\). We do not claim that the
common-target condition holds for every collection of analyses encountered in
practice. Its role is to isolate the effect of search; allowing candidates with
different asymptotic targets would combine that effect with ordinary
candidate-specific misspecification or bias.

\paragraph{Independent evaluation.}
The examples studied in the present paper use the same sampling fluctuations
both to select an analysis and to form the reported conclusion. To see why this
matters, suppose instead that \(D_n\) is used to select \(A_n\), while the selected
analysis is evaluated on an independent sample \(D_n'\). The fluctuations in
\(D_n'\) then play no role in the selection of \(A_n\). Under the common-target
benchmark, this restores target consistency. Indeed, the benchmark implies
\[
\sup_{a\in\mathcal A_n}
\mathbb P\bigl(|T_n(D_n,a)-\theta|>\varepsilon\bigr)\to0
\]
for every \(\varepsilon>0\) (if the supremum did not converge to zero, a
subsequence and a deterministic choice nearly attaining the supremum would
contradict the benchmark). Since \(D_n'\) is independent of \(A_n\) and has the
same distribution as \(D_n\),
\[
\begin{aligned}
\mathbb P\!\left(|T_n(D_n',A_n)-\theta|>\varepsilon\right)
&=
\mathbb E\!\left[
\mathbb P\!\left(
|T_n(D_n',A_n)-\theta|>\varepsilon
\,\middle|\, A_n
\right)
\right] \\
&\leqslant
\sup_{a\in\mathcal A_n}
\mathbb P\!\left(|T_n(D_n',a)-\theta|>\varepsilon\right)
\longrightarrow0.
\end{aligned}
\]
A genuinely fresh holdout thus breaks the link between the fluctuations that
drove selection and those entering the reported conclusion. If the holdout is repeatedly queried and
subsequent analyses are chosen in response to its
results, this independence can be lost, creating
a separate problem of adaptive overfitting
\citep{dwork2015b}. The contrast is therefore not between fixed and
data-dependent selection: \(A_n\) remains data-dependent in both cases. It is
between reporting a conclusion affected by the fluctuations that drove
selection and evaluating the selected analysis on independent data.

\paragraph{Stable pseudotruths.}
A particularly stark separation occurs when
\[
T_n^\star\xrightarrow{P}\theta+b,\qquad b\neq0.
\]
We call \(\theta+b\) a \emph{stable pseudotruth}. Target consistency then fails, while conclusion replicability holds automatically: two independent repetitions converge to the same displaced value, even though the selected analytical paths need not agree and may become increasingly unlikely to coincide. The remainder of the paper identifies regimes
in which search produces such a stable displacement and relates its size to the
geometry of the candidate family.

\section{A minimal Gaussian model: replicable selection away from the target}
\label{sec:minimal}

We now turn to the simplest model in which the distinction introduced in
Section~\ref{sec:framework} can be seen explicitly. The model is deliberately
stripped down: many candidate analyses are available, each is statistically
well behaved when chosen in advance, but the reported analysis is selected
after their outcomes on the observed data have been inspected. Its purpose is
to isolate this search effect before introducing dependence among candidate
analyses.

\subsection{Prespecified and selected analyses}
\label{sec:minimal-selection}

For each \(n\), we specialize the general framework to a simple Gaussian model.
Let the observed data be an \(n\times L_n\) matrix
\[
D_n=(Y_{ij})_{1\leqslant i\leqslant n,\;1\leqslant j\leqslant L_n},
\qquad
Y_{ij}\stackrel{\mathrm{ind}}{\sim}N(\theta,\sigma^2).
\]
Formally, the varying number of columns is compatible with
the fixed-\(P\) framework of Section~\ref{sec:framework}:
one may take each observation to be an infinite sequence of independent
Gaussian coordinates and retain only its first \(L_n\) coordinates. For candidate \(a_{n,j}\), we take \(T_n\) to be the empirical mean of the
\(j\)th column,
\[
T_n(D_n,a_{n,j})
=
\frac{1}{n}\sum_{i=1}^n Y_{ij}
=
\theta+\frac{\sigma}{\sqrt n}Z_{n,j},
\qquad
Z_{n,1},\ldots,Z_{n,L_n}\stackrel{\mathrm{iid}}{\sim}N(0,1).
\]
Every deterministic choice \(j_n\) satisfies
\(T_n(D_n,a_{n,j_n})\xrightarrow{P}\theta\). This verifies the common-target
benchmark of Section~\ref{sec:framework}. The size of \(\mathcal A_n\) is
therefore not by itself a problem; what matters is that the index may be chosen
after the candidate statistics have been observed. We emphasize that this independent Gaussian
model uses \(nL_n\) independent observations and is intended only as a transparent
benchmark for selection among independent candidate statistics, not as a
realistic model of search over a fixed dataset. Section~\ref{sec:correlated}
replaces independence by a shared Gaussian fluctuation structure.

The search rule selects
\[
J_n^\star
=
\operatorname*{arg\,max}_{1\leqslant j\leqslant L_n}T_n(D_n,a_{n,j}),
\qquad
A_n^\star=a_{n,J_n^\star},
\]
and reports
\[
T_n^\star
=
\theta+\frac{\sigma}{\sqrt n}M_{L_n},
\qquad
M_{L_n}:=\max_{1\leqslant j\leqslant L_n}Z_{n,j}.
\]
(Since the candidate statistics are independent and continuously distributed, the maximizing index is almost surely unique.) The candidate laws are unchanged; only the index is now chosen from the realized fluctuations. Any displacement that emerges must therefore
arise from searching over the family rather than
from an intrinsically biased or inconsistent
individual estimator.

\subsection{Exponential search and stable displacement}
\label{sec:exponential-search}

For \(M_L=\max_{1\leqslant j\leqslant L}Z_j\) with independent standard
Gaussians,
\begin{equation}
\label{eq:8}
\frac{M_L}{\sqrt{2\log L}}\xrightarrow{P}1
\qquad(L\to\infty),
\end{equation}
see \citet{leadbetter1983}. The next proposition distinguishes this displacement from the replicability of the selected conclusion.

\begin{proposition}[Exponential analytic search]
\label{prop:exponential-search}
For any sequence of integers \(L_n\geqslant1\), two independent repetitions satisfy
\begin{equation}
\label{eq:10}
T_n^{\star(1)}-T_n^{\star(2)}\xrightarrow{P}0,
\end{equation}
and the corresponding maximizing indices satisfy
\begin{equation}
\label{eq:11}
\mathbb P\!\left(J_n^{\star(1)}=J_n^{\star(2)}\right)=\frac1{L_n}.
\end{equation}
If, in addition, \(\log L_n/n\to c\) for some \(c\in[0,\infty)\), then
\begin{equation}
\label{eq:9}
T_n^\star\xrightarrow{P}\theta+\sigma\sqrt{2c}.
\end{equation}
\end{proposition}

When \(\log L_n=o(n)\), the selected conclusion still converges to \(\theta\). The proof is in Section~\ref{app:proof-exponential}~of the Supplementary Material. In the critical exponential regime \(\log L_n\sim cn\), with \(c>0\), the search converts the individual \(n^{-1/2}\) fluctuations into an order-one displacement.

The replication statement has a complementary interpretation.
Two independent repetitions yield asymptotically close reported
values, even if the selected statistic has no fixed limit.
Yet, when \(L_n\to\infty\), the probability that they select
the same candidate is \(1/L_n\to0\). In the critical regime \(\log L_n/n\to c>0\),
the contrast is particularly striking: both reports
converge to the same displaced limit
\(\theta+\sigma\sqrt{2c}\). Repeating the search rule thus reproduces the displacement induced by selection without
reproducing the selected analysis. This illustrates the distinction between conclusion
replicability, which concerns the reported value, and
path replicability, which concerns the identity of the
selected candidate.

The magnitude of the selection-induced displacement is consistent with the information-theoretic
selection bound of \citet{russo2016}. Here \(I\) denotes mutual information and
\(H\) entropy. Since \(J_n^\star\) is a deterministic function of the data and,
by symmetry, is uniform on \(\{1,\ldots,L_n\}\),
\[
I(J_n^\star;D_n)=H(J_n^\star)=\log L_n.
\]
The sub-Gaussian bound of \citet{russo2016} therefore has scale
\(\sigma\sqrt{2\log L_n/n}\), matching the first-order size of the selection
displacement in this benchmark.

\subsection{Independent search can make the selected conclusion exceptionally stable}
\label{sec:superstable}

The selected maximum not only shifts upward; in the independent model its
fluctuations shrink unusually fast. Indeed, classical Gaussian extreme-value theory
gives, as $L\to\infty$,
\begin{equation}
\label{eq:12}
\operatorname{Var}(M_L)
\sim
\frac{\pi^2}{12\log L},
\end{equation}
and hence
\[
\operatorname{Var}(T_n^\star)
\sim
\frac{\pi^2\sigma^2}{12n\log L_n}.
\]
For two independent repetitions,
\[
\E\!\left[
\bigl(T_n^{\star(1)}-T_n^{\star(2)}\bigr)^2
\right]
=
2\operatorname{Var}(T_n^\star).
\]
The derivation of \eqref{eq:12} is recalled in Section~\ref{app:gaussian-maximum-variance}~of the Supplementary Material.

Consequently, while a prespecified
candidate fluctuates on the $n^{-1/2}$ scale
regardless of $L_n$, the selected conclusion
fluctuates around its own center on the much
smaller $n^{-1}$ scale in the critical regime
$\log L_n\sim cn$, with $c>0$. This distinction matters for replication: numerical agreement across independent studies can become stronger precisely in a regime where the reported conclusion remains displaced from the target. If one observes only the reported numerical conclusions, the search procedure can therefore appear more reproducible than a prespecified analysis. In this example, high numerical replicability is not evidence of greater accuracy; it is a consequence of the selection geometry. The extreme-value facts themselves are classical; our point is their implication for the replication law.

There is a related but distinct connection with \citet{chatterjee2014}, where
superconcentration of Gaussian maxima is linked to chaos of the maximizing
location under small perturbations. In the present paper, path instability is defined across independent repetitions, rather than through local perturbations of a single realization.

\begin{figure}[h]
\centering
\includegraphics[width=0.98\textwidth]{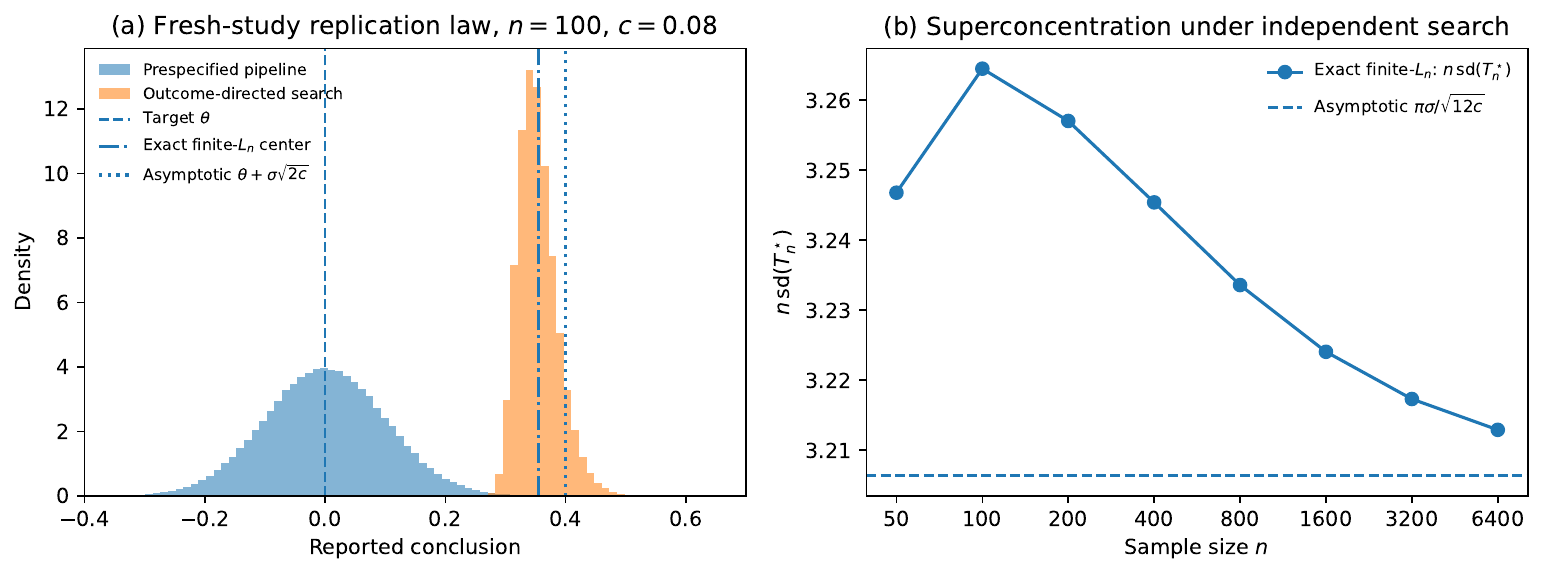}
\caption{Replication behavior in the independent Gaussian model, with
\(L_n\) taken as the nearest integer to \(e^{cn}\). Panel (a) uses
\(300{,}000\) independent repetitions at \(n=100\), with \(\theta=0\),
\(\sigma=1\), and \(c=0.08\), so \(L_n=2981\). The finite-sample expectation
\(\E[T_n^\star]=\sigma\E[M_{L_n}]/\sqrt n\approx0.354\), computed from the exact
distribution of \(M_{L_n}\), is shown alongside the asymptotic limit
\(\sigma\sqrt{2c}=0.400\). Panel (b) shows the finite-sample value of
\(n\operatorname{sd}(T_n^\star)\), obtained from the same distribution, and its
limiting value \(\pi\sigma/\sqrt{12c}\).}
\label{fig:replication-simulation}
\end{figure}

\paragraph{A numerical illustration.}
Figure~\ref{fig:replication-simulation} sets \(\theta=0\), \(\sigma=1\), and
\(c=0.08\), with \(L_n\) the nearest integer to \(e^{cn}\). The histograms in panel (a) are based on simulated
repetitions. The finite-sample means and standard
deviations reported in the text or displayed
in the figure are computed from the Gaussian
model, using in particular the exact distribution
$\mathbb P(M_L\leqslant x)=\Phi(x)^L$. Panel (a) compares a prespecified candidate with the selected report at
\(n=100\); the finite-sample mean of the latter is about \(0.354\), compared
with the asymptotic displacement \(\sqrt{2c}=0.400\). For the same $n=100$ and $L_n=2981$,
the prespecified estimator has an exact finite-sample standard deviation of $0.100$, compared with
approximately $0.033$ for the selected estimator (also computed from its exact distribution). Yet two independent repetitions select the same
candidate with probability only \(1/2981\approx0.00034\). Panel (b) shows the
corresponding finite-sample value of \(n\operatorname{sd}(T_n^\star)\) and its
limit \(\pi/\sqrt{12c}\). Because these quantities are computed from exact formulas, the slight increase between the first two sample sizes is genuine finite-sample behavior, not simulation noise. The
remaining gaps illustrate the relatively slow convergence of Gaussian
extreme-value approximations, even when the asymptotic mechanism is already
visible.

\paragraph{Directional versus two-sided selection.}
The stable displacement above relies on directional selection: the analyst
always favors large positive fluctuations. It is useful to distinguish this
from selection of the most extreme result regardless of sign. If the analyst
instead selects
\[
J_n^{\pm}
=
\operatorname*{arg\,max}_{1\leqslant j\leqslant L_n}
|Z_{n,j}|
\]
and reports
\[
T_n^{\pm}=\theta+\frac{\sigma}{\sqrt n}Z_{n,J_n^{\pm}},
\]
then under \(\log L_n/n\to c>0\),
\[
|T_n^{\pm}-\theta|\xrightarrow{P}\sigma\sqrt{2c},
\]
but symmetry gives
\[
T_n^{\pm}\Rightarrow
\tfrac12\delta_{\theta+\sigma\sqrt{2c}}
+\tfrac12\delta_{\theta-\sigma\sqrt{2c}}.
\]
Thus two-sided selection does not produce the stable signed displacement obtained under directional selection. The same distinction arises when selection is based on statistical significance. For a standardized Gaussian score, the upper-tail
one-sided \(p\)-value is a decreasing function of the signed score, while the
two-sided \(p\)-value is a decreasing function of its absolute value. Hence
selecting the smallest one-sided \(p\)-value is equivalent to maximizing the
signed score, whereas selecting the smallest two-sided \(p\)-value is equivalent
to maximizing its absolute value. Under the centered null \(\theta=0\), the
latter produces the symmetric limiting mixture rather than a single stable
pseudotruth. For two independent repetitions, the selected signs agree with
probability \(1/2\) for every \(n\), so conclusion replicability fails even
though the magnitude is asymptotically replicable.

\paragraph{Stopping at a threshold.}
Maximization provides a transparent benchmark for the selection
mechanism studied in this paper, rather than defining its scope. The same phenomenon can arise, for example, when an analyst examines candidates sequentially and stops at the first one exceeding a fixed threshold. In the independent Gaussian model, fix a reporting threshold \(t\), without requiring
knowledge of \(\theta\), write \(T_{n,j}=T_n(D_n,a_{n,j})\), inspect candidates in a
predetermined order, and report the first \(T_{n,j}\geqslant t\),
using a fixed fallback if no candidate crosses. If \(t>\theta\) and \(\log L_n/n\to
c>(t-\theta)^2/(2\sigma^2)\), a crossing occurs with probability tending to one, and the selected
estimator converges in probability to \(t\). The probability of selecting
the same index in two independent repetitions nevertheless tends to zero. A proof is given in Section~\ref{app:threshold}~of the Supplementary Material.
This idealized rule resembles stopping at the first specification that
meets a substantive reporting threshold.

\subsection{Selection on a proxy criterion}
\label{sec:proxy-selection}

Direct maximization assumes that the selection criterion is the quantity
ultimately reported. In many applications, however, these are different objects: a score may measure fit or
predictive performance, while the quantity ultimately reported estimates the
scientific target. To separate the two, let
\(S_{n,1},\ldots,S_{n,L_n}\stackrel{\mathrm{iid}}{\sim}N(0,1)\) and,
independently, \(\eta_{n,1},\ldots,\eta_{n,L_n}\stackrel{\mathrm{iid}}{\sim}N(0,1)\),
and set
\[
T_{n,j}
=
\theta+\frac{\sigma}{\sqrt n}
\left(\rho_n S_{n,j}+\sqrt{1-\rho_n^2}\,\eta_{n,j}\right),
\qquad |\rho_n|\leqslant1.
\]
The analyst selects \(J_n^\star
=\operatorname*{arg\,max}_{1\leqslant j\leqslant L_n} S_{n,j}\) and reports
\(T_n^\star=T_{n,J_n^\star}\). The parameter \(\rho_n\) is the correlation
between the selection score and the standardized sampling error of the
reported quantity. The cases \(\rho_n=1\) and \(\rho_n=0\) respectively recover direct maximization and selection based on a score independent of the reporting error. Intermediate values describe proxy selection, where the search criterion is
imperfectly aligned with the scientific quantity ultimately reported.

\begin{proposition}[Proxy-directed search]
\label{prop:proxy-search}
Let \(M_{L_n}=\max_{1\leqslant j\leqslant L_n}S_{n,j}\). Then
\begin{equation}
\label{eq:proxy-decomposition}
T_n^\star
\stackrel{d}{=}
\theta+\frac{\sigma}{\sqrt n}
\left(\rho_n M_{L_n}+\sqrt{1-\rho_n^2}\,E_n\right),
\end{equation}
where \(E_n\sim N(0,1)\) is independent of \(M_{L_n}\), and
\begin{equation}
\label{eq:proxy-variance}
\operatorname{Var}(T_n^\star)
=
\frac{\sigma^2}{n}
\left[\rho_n^2\operatorname{Var}(M_{L_n})+1-\rho_n^2\right].
\end{equation}
If \(L_n\to\infty\) and
\begin{equation}
\label{eq:proxy-balance}
\rho_n\sqrt{\frac{2\log L_n}{n}}\longrightarrow\beta\in\mathbb R,
\end{equation}
then
\begin{equation}
\label{eq:proxy-limit}
T_n^\star\xrightarrow{P}\theta+\sigma\beta.
\end{equation}
\end{proposition}

The proof is in Section~\ref{app:proof-proxy}
of the Supplementary Material.
The balance \eqref{eq:proxy-balance}
makes explicit the trade-off between the strength
of the alignment and the size of the search family. When \(\rho_n\to\rho\neq0\) and
\(\log L_n/n\to c>0\), the limiting displacement
is \(\sigma\beta=\sigma\rho\sqrt{2c}\). More surprisingly, a nonzero displacement
remains possible when \(\rho_n\to0\).
Although the selection score and the standardized
reporting error of each fixed candidate then
become asymptotically independent, searching over
a sufficiently large family can amplify this
vanishing dependence. For such a nonzero limit,
\eqref{eq:proxy-balance} requires
\(\rho_n^2\log L_n/n\to\beta^2/2>0\),
and hence \(\log L_n/n\to\infty\). This extreme growth requirement is a limitation
of the independent-candidate model and motivates
the correlated geometry developed next.

The variance formula also shows that the exceptional concentration of direct
maximization is not robust. If \(\rho_n\to\rho\) with \( |\rho|<1\) and
\(L_n\to\infty\), then
\(n\operatorname{Var}(T_n^\star)/\sigma^2\to1-\rho^2\), because
\(\operatorname{Var}(M_{L_n})\to0\). Consequently, in this regime,
the selected conclusion retains the usual
$n^{-1/2}$ standard-deviation scale.
If instead $|\rho_n|\to1$, its concentration
rate depends on the relative sizes of
$1-\rho_n^2$ and $\operatorname{Var}(M_{L_n})$. At the same time, \eqref{eq:proxy-variance} shows that variance need not reveal
the displacement. If \(\rho_n\to0\) while
\(\rho_n\sqrt{2\log L_n/n}\to\beta\neq0\), then
\[
T_n^\star\xrightarrow{P}\theta+\sigma\beta,
\qquad
\operatorname{Var}(T_n^\star)\sim\frac{\sigma^2}{n}.
\]
Thus an order-one search displacement can coexist with the same first-order
variance as a prespecified candidate.

\section{Correlated Gaussian analytic search}
\label{sec:correlated}

The independent model of Section~\ref{sec:minimal} makes the mechanism
transparent, but its independence assumption is much stronger than what one
would expect when many candidate analyses use the same data. Two regression
specifications that differ by one covariate, two nearby transformations of the
same response, or two models fitted to the same variables will typically have
strongly correlated sampling errors. Accordingly, we replace independence by a shared
Gaussian fluctuation structure in which dependence across candidate analyses
is represented geometrically.

This extension concerns a different form of dependence from the proxy model
of Section~\ref{sec:proxy-selection}. There, the parameter \(\rho_n\) measured
the alignment between the fluctuation driving selection and the fluctuation
entering the reported quantity. By contrast, under direct maximization, selection and reporting are driven by the same candidate fluctuation; the question becomes how these fluctuations depend across candidates.

\subsection{Gaussian representation}
\label{sec:gaussian-representation}

Let
\[
G_n\sim N(0,I_{p_n}),
\]
and let \(V_n\subseteq B_2^{p_n}\) be nonempty and compact, where \(B_2^{p_n}\) denotes the closed Euclidean unit ball in \(\mathbb R^{p_n}\). Candidate \(v\in V_n\) reports
\[
T_{n,v}
=
\theta+\frac{\sigma}{\sqrt n}\langle G_n,v\rangle.
\]
Then \(\E T_{n,v}=\theta\),
\(\operatorname{Var}(T_{n,v})=\sigma^2\|v\|_2^2/n\leqslant\sigma^2/n\), and
\[
\operatorname{Cov}(T_{n,u},T_{n,v})
=
\frac{\sigma^2}{n}\langle u,v\rangle.
\]
The geometry of \(V_n\) thus encodes the dependence of the candidate fluctuations. The uniform variance bound also implies
that, for any deterministic sequence $v_n\in V_n$,
\(
T_{n,v_n}\xrightarrow{P}\theta,
\)
so the common-target benchmark remains satisfied throughout the correlated
model. Here \(p_n\) is the dimension of the underlying Gaussian
fluctuation space, not the number of candidate analyses. The latter may be
finite or infinite. In particular, the independent benchmark with \(L_n\) candidates is recovered by taking \(p_n=L_n\) and
\(
V_n=\{e_1,\ldots,e_{L_n}\},
\)
where \(e_1,\ldots,e_{L_n}\) are the canonical basis vectors of \(\mathbb R^{L_n}\). Nearby directions have similar sampling fluctuations, whereas nearly
orthogonal directions have weakly correlated fluctuations. In this way the
inner-product geometry of \(V_n\) records which analyses offer genuinely
different opportunities to exploit sampling noise. 

The covariance matrix $\Sigma$ of any
finite centered Gaussian vector with marginal
variances at most one is positive semidefinite
and can therefore be written as a Gram matrix
$\Sigma_{jk}=\langle v_j,v_k\rangle$,
with $\|v_j\|_2\leqslant1$.
Thus the Gaussian representation introduced above
entails no loss of generality for such vectors. Beyond the finite-dimensional setting, the model may be read as a Gaussian
approximation to a growing family of estimators. In many classical problems a
fixed collection of regular estimators has a joint Gaussian limit after
\(\sqrt n\) rescaling. Here the collection itself is allowed to grow with
\(n\), so a concrete application requires the Gaussian approximation to remain
valid uniformly over the candidate family. For maxima of high-dimensional sums, see \citet{chernozhukov2013}.

Under direct maximization, the search rule reports
\begin{equation}
\label{eq:17}
T_n^\star
=
\sup_{v\in V_n}T_{n,v}
=
\theta+\frac{\sigma}{\sqrt n}
\sup_{v\in V_n}\langle G_n,v\rangle.
\end{equation}
Compactness ensures that the supremum is attained. The analyst need not observe
\(G_n\) itself; it is only a representation of the joint sampling law of the
candidate statistics. What the analyst observes are the candidate results
computed from the data, and those results determine the selected analysis. The key difference from the independent model is that search complexity can no
longer be summarized by the number of candidate analyses alone. If many candidates are nearly identical, adding more of them creates little additional scope for exploiting sampling noise; candidates that are more widely separated in the Gaussian geometry can contribute substantially more.

\paragraph{A regression interpretation.}
The Gaussian representation has a direct interpretation in marginal variable
screening. Under the global null, let \(Y\sim N(0,\sigma^2 I_n)\) and
consider fixed regressors \(x_1,\ldots,x_{L_n}\in\mathbb R^n\) with
\(\|x_j\|_2^2=n\). For each \(j\), fit a separate univariate
least-squares regression of \(Y\) on \(x_j\), without an intercept.
In the notation above, take \(p_n=n\), \(G_n=Y/\sigma\), \(\theta=0\),
and \(V_n=\{v_1,\ldots,v_{L_n}\}\), where \(v_j=x_j/\sqrt n\).
Then \(T_{n,v_j}=x_j^\top Y/n=\widehat\beta_j\) is the marginal
least-squares estimator, and
\(\operatorname{Cov}(T_{n,v_j},T_{n,v_k})
=\sigma^2 x_j^\top x_k/n^2\).
Thus maximizing \(T_{n,v_j}\) selects the largest signed marginal
coefficient, with dependence determined by the regressors' Gram matrix.
The same Gaussian representation also covers selection by the largest absolute marginal coefficient, as when identifying the first variable to enter a no-intercept Lasso path with standardized regressors. This corresponds to maximizing over \(V_n\cup(-V_n)\), rather than \(V_n\), and the maximum is the magnitude of the selected coefficient. Reporting the selected coefficient with its
original sign instead leads to a two-sided
selection-and-reporting rule, analogous to the one studied in
Section~\ref{sec:superstable} for independent
Gaussian candidates. The example shows that the geometry underlying our replication results is
not merely an abstract construction: it already appears in the familiar
practice of screening covariates one at a time and choosing which result
to report.

\subsection{Gaussian width}
\label{sec:gaussian-width}

Define the Gaussian width
\begin{equation}
\label{eq:18}
w(V_n)
=
\E\left[\sup_{v\in V_n}\langle G_n,v\rangle\right],
\end{equation}
and write \(w_n=w(V_n)\). Gaussian width is a classical measure of effective
size that automatically accounts for correlation
\citep{ledoux1991,vershynin2018}. For a realization \(G_n\), the
supremum in \eqref{eq:18} is the largest normalized sampling fluctuation
available to the analyst; its expectation over repeated datasets is exactly
\(w_n\). In the present setting, this gives the exact mean displacement:
\[
\E T_n^\star=\theta+\frac{\sigma w_n}{\sqrt n}.
\]
The main question is therefore how tightly the selected statistic concentrates
around this displaced center.

\begin{theorem}[Correlated Gaussian search and replication]
\label{thm:gaussian-search}
Let \(V_n\subseteq B_2^{p_n}\) be nonempty and compact, and let \(T_n^\star\)
be defined by \eqref{eq:17}. Then, for every \(t>0\),
\begin{equation}
\label{eq:20}
\mathbb P\!\left(
\left|T_n^\star-\theta-\frac{\sigma w_n}{\sqrt n}\right|>t
\right)
\leqslant
2\exp\!\left(-\frac{nt^2}{2\sigma^2}\right).
\end{equation}
Moreover, if \(T_n^{\star(1)}\) and \(T_n^{\star(2)}\) are two independent
repetitions, then
\begin{equation}
\label{eq:21}
\E\!\left[\left(T_n^{\star(1)}-T_n^{\star(2)}\right)^2\right]
\leqslant\frac{2\sigma^2}{n}.
\end{equation}
In particular,
\[
T_n^{\star(1)}-T_n^{\star(2)}\xrightarrow{P}0.
\]
If, in addition,
\begin{equation}
\label{eq:22}
\frac{w_n^2}{n}\longrightarrow\alpha\in[0,\infty),
\end{equation}
then
\begin{equation}
\label{eq:23}
T_n^\star\xrightarrow{P}\theta+\sigma\sqrt\alpha.
\end{equation}
\end{theorem}

\begin{proof}
Let \(F_n(g)=\sup_{v\in V_n}\langle g,v\rangle\). Since
\(V_n\subseteq B_2^{p_n}\),
\[
|F_n(g)-F_n(g')|
\leqslant
\sup_{v\in V_n}|\langle g-g',v\rangle|
\leqslant
\|g-g'\|_2,
\]
so \(F_n\) is \(1\)-Lipschitz. The Gaussian concentration inequality for
Lipschitz functions \citep{ledoux1991} therefore gives
\[
\mathbb P\!\left(|F_n(G_n)-\E F_n(G_n)|>s\right)
\leqslant2e^{-s^2/2},
\]
which yields \eqref{eq:20} after taking \(s=\sqrt n\,t/\sigma\). The Gaussian Poincar\'e inequality \citep{boucheron2013} gives
\(\operatorname{Var}(F_n(G_n))\leqslant1\), so
\(\operatorname{Var}(T_n^\star)\leqslant\sigma^2/n\).
For two independent copies,
\[
\mathbb E\!\left[\left(T_n^{\star(1)}-T_n^{\star(2)}\right)^2\right]
=
2\operatorname{Var}(T_n^\star),
\]
which gives \eqref{eq:21} and, in particular,
\[
T_n^{\star(1)}-T_n^{\star(2)}\xrightarrow{P}0.
\]
Finally, \(w_n\geqslant0\), so \eqref{eq:22} implies
\(w_n/\sqrt n\to\sqrt\alpha\); combining this with \eqref{eq:20} proves
\eqref{eq:23}.
\end{proof}

Theorem~\ref{thm:gaussian-search} is the central scale statement of the paper:
\[
T_n^\star
=
\theta+\frac{\sigma w_n}{\sqrt n}+O_P(n^{-1/2}).
\]
For the independent family
$V_n=\{e_1,\ldots,e_{L_n}\}$, we have $w_n=\E M_{L_n}$, with
$0\leqslant w_n\leqslant \sqrt{2\log L_n}$ for every $L_n\geqslant 1$ and,
by classical Gaussian extreme-value theory,
$w_n\sim\sqrt{2\log L_n}$ when $L_n\to\infty$.
Consequently, if $\log L_n/n\to c\in[0,\infty)$,
then $w_n/\sqrt n\to\sqrt{2c}$,
recovering the limiting displacement in
Proposition~\ref{prop:exponential-search}.
The theorem also recovers the proposition's unconditional
conclusion replicability. More generally,
\(w_n/\sqrt n\) determines the scale of the search displacement, and
\(w_n^2\) of order \(n\) is the critical scale for an order-one displacement.
In particular, if \(w_n^2/n\to\alpha>0\), then the replication law converges
to a point mass at \(\theta+\sigma\sqrt\alpha\). The complete procedure
therefore has a stable pseudotruth even though every deterministic sequence
of candidates is consistent for \(\theta\).

The variance bound in Theorem~\ref{thm:gaussian-search} is only an upper
bound. The geometry of \(V_n\) determines whether the selected conclusion is
superconcentrated or retains the ordinary \(n^{-1/2}\) fluctuation scale, and
there is no matching universal lower bound. For independent search with \(\log L_n\) of order \(n\),
the variance can be of order \(n^{-2}\), whereas the hypercube and spherical
examples of Section~\ref{sec:search-dimension} retain variance of order
\(n^{-1}\). Thus stable displacement and conclusion replicability do not
require superconcentration.

Since \(\sqrt n(T_n^\star-\theta)/\sigma=w_n+O_P(1)\), if \(w_n\to\infty\),
any naive interval \(T_n^\star\pm z\sigma/\sqrt n\), with fixed \(z>0\),
excludes \(\theta\) with probability tending to one. By contrast, the theorem
does not imply any general conclusion about path agreement: whether independent
repetitions select the same candidate depends on the geometry of \(V_n\), a
question we return to in Section~\ref{sec:search-dimension}.

Finally, the maximizing rule also provides an envelope for any selection rule taking values in \(V_n\): pointwise,
\(T_{n,A_n}\leqslant T_n^\star\), and hence
\(\E T_{n,A_n}-\theta\leqslant\sigma w_n/\sqrt n\). Thus Gaussian width bounds the largest possible upward mean displacement within the candidate family.

\section{Search dimension and metric geometry}
\label{sec:search-dimension}

Theorem~\ref{thm:gaussian-search} suggests the effective dimension
\[
d_{\rm search}(V):=w(V)^2.
\]
This is simply the squared Gaussian width; it need not be an integer or
coincide with the ambient dimension. Since
\(
\E T_n^\star-\theta=\sigma\sqrt{d_{\rm search}(V_n)/n},
\)
\(d_{\rm search}(V_n)\) is the squared expected selection displacement measured
in units of the reference variance scale \(\sigma^2/n\). To see how dependence
changes this effective dimension, we compute it for several candidate
geometries.

\subsection{Equicorrelated pipelines}

The simplest way to see how correlation changes analytical search is to keep
the marginal variance of every candidate unchanged while introducing a common
component into their sampling errors. Let
\[
X_{n,j}=\sqrt{\rho_n}\,Z_0+\sqrt{1-\rho_n}\,Z_j,
\qquad j=1,\ldots,L_n,
\]
where \(Z_0,Z_1,\ldots,Z_{L_n}\) are independent standard Gaussians and
\(0\leqslant\rho_n<1\). Then each candidate has unit variance and pairwise
correlation \(\rho_n\). Candidate \(j\) reports
\(T_{n,j}=\theta+\sigma X_{n,j}/\sqrt n\). This corresponds to
\(p_n=L_n+1\) and
\(V_n=\{v_{n,1},\ldots,v_{n,L_n}\}\subset\mathbb R^{L_n+1}\), with
\(v_{n,j}=\sqrt{\rho_n}\,e_0+\sqrt{1-\rho_n}\,e_j\) and
\(G_n=(Z_0,Z_1,\ldots,Z_{L_n})^\top\). With
\(M_L=\max_{1\leqslant j\leqslant L}Z_j\) and \(m_L=\E M_L\),
\[
\max_{1\leqslant j\leqslant L_n}X_{n,j}
=
\sqrt{\rho_n}\,Z_0+\sqrt{1-\rho_n}\,M_{L_n},
\qquad
w(V_n)=\sqrt{1-\rho_n}\,m_{L_n}.
\]
The decomposition makes the role of correlation explicit. The common
fluctuation moves all candidates together and cannot be amplified by selecting
among them; only the candidate-specific fluctuations contribute to the
search-induced displacement. Hence, if \(L_n\to\infty\),
\[
d_{\rm search}(V_n)
\sim
2(1-\rho_n)\log L_n.
\]
The common fluctuation survives, however, in replication variability:
\[
\E\!\left[
\bigl(T_n^{\star(1)}-T_n^{\star(2)}\bigr)^2
\right]
=
\frac{2\sigma^2}{n}
\left[\rho_n+(1-\rho_n)\operatorname{Var}(M_{L_n})\right].
\]
Thus correlation reduces the effective search scale by the factor
\(1-\rho_n\), and the critical balance becomes
\((1-\rho_n)\log L_n\) of order \(n\). If \(\rho_n\to\rho>0\), the candidate-specific contribution
\(\operatorname{Var}(M_{L_n})\) vanishes while the common component remains,
so \(\operatorname{Var}(T_n^\star)\sim\sigma^2\rho/n\), meaning that the selected conclusion fluctuates on the ordinary
$n^{-1/2}$ scale, rather than the smaller scale of the
independent maximum. The common
component does not affect the maximizing index, which remains uniform, so two
independent repetitions agree on the selected pipeline with probability
\(1/L_n\). Correlation therefore reduces the displacement generated by search
without making the winning path any more stable.

\subsection{A combinatorial analysis space: the hypercube}
\label{sec:hypercube}

Real analytical search spaces may be built from many elementary choices---for
example, whether to include particular variables, transformations, or modeling
components---so that a moderate number of local decisions generates a very
large family of complete pipelines. A simple model of this combinatorial
structure is the hypercube.

For \(p\geqslant1\), consider
\begin{equation}
\label{eq:37}
V_p
=
\left\{
\frac1{\sqrt p}(\varepsilon_1,\ldots,\varepsilon_p)^\top:
\varepsilon_k\in\{-1,+1\}
\right\}.
\end{equation}
The family contains \(2^p\) unit-norm pipelines generated by \(p\) binary
choices. If two sign vectors differ in \(h\) coordinates, their inner product
is \(1-2h/p\), so the correlation varies across the hypercube. Pipelines
differing in only a few decisions are strongly positively correlated, those
at Hamming distance \(p/2\), when \(p\) is even, are orthogonal, and opposite
pipelines have correlation \(-1\). No single correlation parameter summarizes
this geometry. Since each sign may be optimized coordinatewise,
\begin{equation}
\label{eq:38}
\sup_{v\in V_p}\langle G,v\rangle
=
\frac1{\sqrt p}\sum_{k=1}^p|G_k|.
\end{equation}
Thus optimization over \(2^p\) complete pipelines decomposes into \(p\)
coordinatewise choices, with each sign chosen to match the corresponding realized Gaussian fluctuation.

\begin{proposition}[Hypercube search]
\label{prop:hypercube}
For the family \(V_p\) in \eqref{eq:37}, with \(G\sim N(0,I_p)\),
\begin{equation}
\label{eq:39}
w(V_p)=\sqrt{\frac{2p}{\pi}},
\qquad
d_{\rm search}(V_p)=\frac{2p}{\pi},
\end{equation}
and
\begin{equation}
\label{eq:40}
\operatorname{Var}\!\left(\sup_{v\in V_p}\langle G,v\rangle\right)
=
1-\frac{2}{\pi}.
\end{equation}
For any sequence of integers \(p_n\geqslant1\),
two independent repetitions satisfy
\begin{equation}
\label{eq:42}
\E\!\left[\bigl(T_n^{\star(1)}-T_n^{\star(2)}\bigr)^2\right]
=
\frac{2\sigma^2}{n}\left(1-\frac{2}{\pi}\right).
\end{equation}
If \(A_n^\star\) denotes the maximizing pipeline, then
\begin{equation}
\label{eq:43}
\mathbb P\!\left(A_n^{\star(1)}=A_n^{\star(2)}\right)=2^{-p_n}.
\end{equation}
If, in addition,
\(p_n/n\to\kappa\in[0,\infty)\), then
\begin{equation}
\label{eq:41}
T_n^\star\xrightarrow{P}
\theta+\sigma\sqrt{\frac{2\kappa}{\pi}}.
\end{equation}
\end{proposition}

The proof is in Section~\ref{app:hypercube-proof}~of the Supplementary Material. The example shows that the
relevant complexity is linear in the number \(p\) of elementary binary choices,
and hence logarithmic in the number \(2^p\) of complete pipelines:
\[
d_{\rm search}(V_p)
=
\frac{2}{\pi\log 2}\log |V_p|.
\]
Unlike the independent exponential-search model, the selected value retains
the ordinary \(n^{-1/2}\) fluctuation scale. When \(p_n/n\to\kappa>0\), both studies
converge to the same displaced limit, while the probability of choosing the
same vertex goes to zero. The optimized value concentrates because it is a
normalized sum of many independent coordinatewise contributions, even though the signs defining the selected pipeline typically differ across studies. Note that the same path instability occurs when \(p_n\to\infty\) but
\(p_n/n\to0\), in which case both reported conclusions converge to \(\theta\).

\subsection{A continuum of analyses}
\label{sec:sphere}

The phenomenon does not require a discrete search space. Consider the unit sphere \(S^{p-1}\subset\mathbb R^p\) for \(p\geqslant2\), so the analyst may choose among a continuum of candidate directions. Every fixed direction still has variance \(\sigma^2/n\) and is
consistent for \(\theta\); selection instead aligns the chosen direction with
the realized high-dimensional sampling fluctuation.

By Cauchy--Schwarz,
\[
\sup_{v\in S^{p-1}}\langle G,v\rangle=\|G\|_2,
\qquad
A^\star=\frac{G}{\|G\|_2}.
\]
Hence, as \(p\to\infty\),
\[
w(S^{p-1})
=
\E\|G\|_2
=
\sqrt2\,\frac{\Gamma((p+1)/2)}{\Gamma(p/2)}
\sim\sqrt p,
\qquad
d_{\rm search}(S^{p-1})\sim p.
\]
If \(p_n/n\to\kappa\in[0,\infty)\),  then
\[
T_n^\star\xrightarrow{P}\theta+\sigma\sqrt\kappa.
\]
Thus the critical scale is again \(p_n\) of order \(n\): below this scale the
search displacement vanishes, while at this scale it can remain of constant
order. In the critical regime \(\kappa>0\), we necessarily have \(p_n\to\infty\),
and
\(\operatorname{Var}(\|G_n\|_2)\to1/2\), so
\(\operatorname{Var}(T_n^\star)\sim\sigma^2/(2n)\). Unlike the independent
maximum, the selected conclusion therefore retains the ordinary
\(n^{-1/2}\) standard-deviation scale. The short derivations are collected in
Section~\ref{app:sphere-calculations}~of the Supplementary Material.

The selected paths in two independent repetitions are independent uniform
directions on \(S^{p_n-1}\). For \(p_n\geqslant2\), exact equality has
probability zero, and
\[
\E\langle A_n^{\star(1)},A_n^{\star(2)}\rangle=0,
\qquad
\E\langle A_n^{\star(1)},A_n^{\star(2)}\rangle^2=\frac1{p_n}.
\]
Thus the selected directions become nearly orthogonal when \(p_n\to\infty\),
while their optimized norms satisfy
\(\|G_n\|_2/\sqrt{p_n}\xrightarrow{P}1\). Geometrically, each study aligns its analysis with its own realization of the sampling noise: the selected directions differ substantially across studies, while the norms of the underlying Gaussian vectors concentrate. The sphere makes especially clear why cardinality is not
the relevant notion of analytical complexity---the candidate family is
uncountable, yet its search dimension is only of order \(p\).

\subsection{Packing, covering, and generic chaining}
\label{sec:chaining}

The examples suggest a formulation that does not rely on Euclidean
coordinates. Let \(\{X_{n,a}:a\in\mathcal A_n\}\) be a centered separable
Gaussian process such that
\(\sup_a\operatorname{Var}(X_{n,a})\leqslant1\) and
\(w_n:=\E\sup_a X_{n,a}<\infty\). For fixed \(\sigma>0\), let candidate
\(a\) report
\[
T_{n,a}
=
\theta+\frac{\sigma}{\sqrt n}X_{n,a},
\]
so every deterministic sequence of candidates satisfies the common-target
benchmark. Under direct maximization,
\[
T_n^\star
=
\theta+\frac{\sigma}{\sqrt n}\sup_{a\in\mathcal A_n}X_{n,a}.
\]
The process has canonical pseudometric
\[
d_n(a,b)
=
\left\{\E(X_{n,a}-X_{n,b})^2\right\}^{1/2}.
\]
This distance describes how differently two analyses react to sampling
noise. In the Euclidean model, \(X_{n,v}=\langle G_n,v\rangle\), it
reduces to \(d_n(u,v)=\|u-v\|_2\). Borell--TIS concentration
\citep{borell1975,tsirelson1976} gives
\[
T_n^\star
=
\theta+\frac{\sigma w_n}{\sqrt n}+O_P(n^{-1/2}).
\]
Thus \(w_n\) measures the mean displacement induced by search, while the
remaining random fluctuation is at most on the \(n^{-1/2}\) scale. In this setting we write \(d_{\rm search}(\mathcal A_n):=w_n^2\) for the corresponding search dimension.

Let \(M(\mathcal A,d,\varepsilon)\) be the largest number of points with
pairwise \(d\)-distance exceeding \(\varepsilon\), and let
\(N(\mathcal A,d,\varepsilon)\) be the smallest number of \(d\)-balls of
radius \(\varepsilon\) needed to cover \(\mathcal A\). These are the packing
and covering numbers; \(\log N\) is the metric entropy at resolution
\(\varepsilon\).

\begin{proposition}[Sudakov--Dudley bounds]
\label{prop:entropy-bounds}
There exist universal constants \(c_0,C_0>0\) such that
\[
c_0\sup_{\varepsilon>0}\varepsilon
\sqrt{\log M(\mathcal A_n,d_n,\varepsilon)}
\leqslant w_n
\leqslant
C_0\int_0^{\operatorname{diam}(\mathcal A_n,d_n)}
\sqrt{\log N(\mathcal A_n,d_n,\varepsilon)}\,d\varepsilon.
\]
\end{proposition}

The bounds are the classical Sudakov minoration and Dudley entropy integral;
see \citet{dudley1967} and \citet{talagrand2014}. Since
\(d_{\rm search}(\mathcal A_n)=w_n^2\), they relate the search dimension
directly to the metric geometry of the candidate family. A large packing at a
given resolution provides many sufficiently separated sampling fluctuations
that search can exploit, whereas small covering numbers across all resolutions
keep the search dimension under control. Exponentially many nominal pipelines
matter only to the extent that their sampling fluctuations remain sufficiently
separated. If many pipelines produce nearly the same fluctuation pattern,
their multiplicity contributes little additional search dimension.

Write
\[
\mathcal J_n:=\int_0^{\operatorname{diam}(\mathcal A_n,d_n)}
\sqrt{\log N(\mathcal A_n,d_n,\varepsilon)}\,d\varepsilon.
\]

\begin{corollary}[Geometric search regimes]
\label{cor:geometric-regimes}
If \(\mathcal J_n=o(\sqrt n)\), then \(T_n^\star\xrightarrow{P}\theta\). Conversely, if
there exist \(\varepsilon_0,\eta>0\) such that
\begin{equation}
\label{eq:60}
\liminf_{n\to\infty}\frac1n
\log M(\mathcal A_n,d_n,\varepsilon_0)\geqslant\eta,
\end{equation}
then
\begin{equation}
\label{eq:61}
\liminf_{n\to\infty}\frac{w_n}{\sqrt n}
\geqslant c_0\varepsilon_0\sqrt\eta,
\end{equation}
and therefore
\begin{equation}
\label{eq:62}
\liminf_{n\to\infty}(\E T_n^\star-\theta)
\geqslant\sigma c_0\varepsilon_0\sqrt\eta.
\end{equation}
Regardless of these regimes, Gaussian concentration gives, for two independent
repetitions,
\[
T_n^{\star(1)}-T_n^{\star(2)}\xrightarrow{P}0.
\]
\end{corollary}

The proof is in Section~\ref{app:entropy-corollary}~of the Supplementary Material. Condition \eqref{eq:60} is sufficient for a nonvanishing displacement but does
not identify its limiting value; that requires control of the normalized width
itself. It is also a one-scale condition. Gaussian complexity may accumulate
across many resolutions even when no fixed scale contains an exponentially
large packing. This multiscale phenomenon leads naturally to generic chaining.

For admissible finite
sets \(\mathcal A_{n,k}\subseteq\mathcal A_n\) satisfying
\(|\mathcal A_{n,0}|=1\) and
\(|\mathcal A_{n,k}|\leqslant2^{2^k}\), define
\[
\gamma_2(\mathcal A_n,d_n)
=
\inf_{\{\mathcal A_{n,k}\}}
\sup_{a\in\mathcal A_n}
\sum_{k\geqslant0}2^{k/2}d_n(a,\mathcal A_{n,k}),
\]
where \(d_n(a,\mathcal A_{n,k})
=\inf_{b\in\mathcal A_{n,k}}d_n(a,b)\).
Rather than measure complexity at one fixed resolution, \(\gamma_2\)
aggregates the difficulty of approximating the analysis space across all
resolutions. Talagrand's majorizing-measure theorem \citep{talagrand2014} states that there exist
universal constants \(0<c<C<\infty\) such that
\[
c\,\gamma_2(\mathcal A_n,d_n)
\leqslant w_n
\leqslant C\,\gamma_2(\mathcal A_n,d_n).
\]
Thus \(d_{\rm search}(\mathcal A_n)=w_n^2\) is, up to universal constants,
of the same order as \(\gamma_2(\mathcal A_n,d_n)^2\), and the intrinsic
critical scale is \(\gamma_2(\mathcal A_n,d_n)^2\) of order \(n\). Since \(w_n^2\sim 2\log L_n\) in the independent benchmark as \(L_n\to\infty\), the relation above shows that \(\gamma_2^2\) plays, for a general Gaussian process, the role played by \(\log L_n\) in the independent finite model. It thus quantifies the effective scope for exploiting sampling noise while remaining meaningful for correlated, continuous, and multiscale candidate families.

\paragraph{Connection with statistical dimension.}
Returning to the Euclidean Gaussian model, search dimension is closely
related to statistical dimension for conic families. If
\(C\subseteq\mathbb R^p\) is a nontrivial closed convex cone, that is, a closed convex set \(C\neq\{0\}\)
such that \(x\in C\) implies \(\lambda x\in C\) for every
\(\lambda\geqslant0\), let
\[
V=C\cap S^{p-1},
\qquad
\delta(C)=\E\|\Pi_C(G)\|_2^2,
\qquad G\sim N(0,I_p).
\]
The quantity \(\delta(C)\) is called
the statistical dimension of \(C\). With \(w(V)=\E\sup_{v\in V}\langle G,v\rangle\),
\[
d_{\rm search}(V)
=
w(C\cap S^{p-1})^2
\leqslant\delta(C)
\leqslant
w(C\cap S^{p-1})^2+1
=
d_{\rm search}(V)+1
\]
\citep{amelunxen2014}. Thus, up to an additive constant, statistical dimension coincides with the search dimension of the corresponding spherical section. For an \(r_n\)-dimensional linear subspace \(C_n\), \(\delta(C_n)=r_n\), so
the corresponding critical scale is \(r_n\) of order \(n\).

\section{Beyond Gaussianity: search entropy and large deviations}
\label{sec:large-deviations}

The Gaussian theory gives a sharp treatment of correlation because Gaussian
width and the canonical metric encode the joint sampling fluctuations exactly.
The basic replication mechanism, however, is not intrinsically Gaussian.
Return to independent and identically distributed candidate statistics
\(T_{n,1},\ldots,T_{n,L_n}\), all satisfying the common-target benchmark, and
suppose that their upper tails satisfy the large-deviation relation
\citep{dembozeitouni1998}
\begin{equation}
\label{eq:70}
-\frac1n\log\mathbb P(T_{n,1}\geqslant t)\longrightarrow I(t),
\qquad t\in(\theta,t_+),
\end{equation}
where \(I\) is continuous on \([\theta,t_+)\), strictly increasing on
\((\theta,t_+)\), and \(I(\theta)=0\). The rate \(I(t)\) measures the
exponential sampling cost of obtaining a prespecified candidate conclusion
at least as large as \(t\). Write \(I(t_+^-)=\lim_{t\uparrow t_+}I(t)\). If
\(\log L_n/n\to c>0\) with
\[
0<c<I(t_+^-),
\]
let \(t_c\in(\theta,t_+)\) be the unique solution of
\begin{equation}
\label{eq:72}
I(t_c)=c.
\end{equation}

\begin{proposition}[Search entropy versus large-deviation cost]
\label{prop:large-deviation-search}
For \(T_n^\star=\max_{1\leqslant j\leqslant L_n}T_{n,j}\),
\begin{equation}
\label{eq:73}
T_n^\star\xrightarrow{P}t_c.
\end{equation}
Consequently, two independent repetitions satisfy
\(T_n^{\star(1)}-T_n^{\star(2)}\xrightarrow{P}0\).
If the common candidate distribution is continuous and
\(J_n^\star=\operatorname*{arg\,max}_{1\leqslant j\leqslant L_n}T_{n,j}\),
then
\[
\mathbb P\!\left(J_n^{\star(1)}=J_n^{\star(2)}\right)=\frac1{L_n}.
\]
\end{proposition}

The proof is in Section~\ref{app:proof-large-deviation}~of the Supplementary Material. The selected level
has a simple interpretation as a balance between search opportunity and
sampling rarity. A level \(t\) costs approximately \(e^{-nI(t)}\) for one
prespecified candidate, while the analyst has approximately \(e^{nc}\)
opportunities to reach it. The expected number of exceedances therefore has exponential rate
\(c-I(t)\): levels with \(I(t)>c\) are too rare to be reached, while levels
with \(I(t)<c\) occur many times. The selected
conclusion concentrates at the boundary \(I(t)=c\). For Gaussian means, this gives
\(I(t)=(t-\theta)^2/(2\sigma^2)\) and recovers
\(\theta+\sigma\sqrt{2c}\).

A concrete non-Gaussian example is given by independent sample means of
Bernoulli observations. Specifically, let
\[
T_{n,j}=\frac1n\sum_{i=1}^n Y_{ij},
\qquad
Y_{ij}\stackrel{\mathrm{ind}}{\sim}\operatorname{Bernoulli}(\theta),
\qquad 0<\theta<1.
\]
Cram\'er's theorem gives, for \(t\in(\theta,1)\),
\[
I_\theta(t)
=
t\log\frac{t}{\theta}+(1-t)\log\frac{1-t}{1-\theta}.
\]
Since \(I_\theta(1)=\log(1/\theta)\), for every
\(0<c<\log(1/\theta)\) the maximum over \(\exp\{cn+o(n)\}\) individually
consistent means converges to the unique \(t_c\) with \(I_\theta(t_c)=c\).
Ties can occur because the means are lattice-valued, so path replication then
depends on the tie-breaking rule, which is
encoded in the selection kernel \(Q_n\).

Thus the stable-pseudotruth mechanism is not a peculiarity of Gaussian maxima.
What is Gaussian-specific is the sharp geometric treatment of dependence, not
the entropy--rarity balance underlying independent search. In this sense,
Gaussian width and the large-deviation entropy--rate balance provide two
versions of the same principle: the richness of the search determines how far
into the favorable tail a repeated procedure can reliably reach.

\section{Discussion and limitations}
\label{sec:scope}

The main lesson of the preceding results is that conclusion replicability
and target consistency answer different questions. In the Gaussian model,
search shifts the center of the reported conclusion by
\(\sigma w_n/\sqrt n\), while its fluctuations around that center are at
most of order \(n^{-1/2}\). Consequently, independent repetitions may
produce nearly identical conclusions even when those conclusions remain
persistently displaced from \(\theta\). This distinction is central to
interpreting what agreement across repeated studies does, and does not,
establish.

The number of candidate analyses is not, by itself, an intrinsic measure
of search complexity. Independent candidates, highly correlated
specifications, combinatorial families, and continua of directions may
have very different cardinalities yet comparable Gaussian widths.
In the Gaussian setting, packing, covering, and generic chaining relate
this complexity to the canonical sampling metric. For independent
non-Gaussian candidates, the analogous balance is between search entropy
and large-deviation cost, with the selected level determined by \(I(t)=c\).

The theory suggests two practical diagnostics. First, conclusion agreement
should be examined jointly with path agreement. Independent studies that
produce similar numerical conclusions while selecting very different
pipelines display the pattern highlighted by the hypercube and spherical
examples. Such a pattern is not, by itself, evidence of a stable pseudotruth:
similar conclusions may also arise when both studies consistently estimate
the scientific target. Assessing path stability also requires a scientifically
meaningful notion of distance between pipelines and, for formal inference,
an appropriate reference distribution. These considerations argue for
recording the selected pipeline, and preferably the set of candidate analyses considered during the search, alongside the final conclusion. This practice is in the spirit of multiverse analysis
\citep{steegen2016} and complements the stability perspective
of \citet{yukumbier2020} by distinguishing stability of the
report from stability of the analytical path.
Related many-analyst studies examine variation across
analyses of the same data \citep{aczel2026}, whereas our
replication experiment uses independent samples.

Second, repeated Gaussian search can produce unusually small between-study
variation, even when the reported conclusions remain displaced from the
scientific target. Consider \(K\geqslant2\) independent repetitions of the
direct-maximization procedure in the correlated Gaussian model of
Section~\ref{sec:correlated}. For each repetition \(k=1,\ldots,K\), let
\[
T_{n,k}^\star
=
\theta+\frac{\sigma}{\sqrt n}
\sup_{v\in V_n}\langle G_n^{(k)},v\rangle,
\]
where \(V_n\subseteq B_2^{p_n}\) and the \(G_n^{(k)}\)'s are independent
standard Gaussian vectors. By Theorem~\ref{thm:gaussian-search}, each
selected report has expectation \(\theta+\sigma w_n/\sqrt n\) and variance
at most \(\sigma^2/n\).

To see what happens when studies are aggregated, write
\[
T_{n,k}^\star
=
\theta+\frac{\sigma w_n}{\sqrt n}+R_{n,k},
\]
where \(\E R_{n,k}=0\) and
\(\operatorname{Var}(R_{n,k})\leqslant\sigma^2/n\).
For
\(\overline T_{n,K}^\star=K^{-1}\sum_{k=1}^K T_{n,k}^\star\),
independence gives
\[
\overline T_{n,K}^\star-\theta
=
\frac{\sigma w_n}{\sqrt n}
+\frac1K\sum_{k=1}^K R_{n,k},
\qquad
\operatorname{Var}(\overline T_{n,K}^\star)
\leqslant\frac{\sigma^2}{nK}.
\]
If each selected report is naively treated as a prespecified unit-norm
Gaussian estimator with standard error \(\sigma/\sqrt n\), the
standardized aggregate satisfies
\[
\frac{\overline T_{n,K}^\star-\theta}{\sigma/\sqrt{nK}}
=
\sqrt K\,w_n+Z_{n,K},
\]
where \(Z_{n,K}\) has mean zero and variance at most one, and is therefore
\(O_P(1)\). In the critical regime \(w_n^2/n\to\alpha>0\), the
standardized displacement diverges, while for every fixed \(K\),
\(\overline T_{n,K}^\star\) still converges to
\(\theta+\sigma\sqrt\alpha\). Averaging thus reduces random uncertainty
around the selected center without removing its displacement. Consequently, a naive test of the true target based on this standardized
statistic would reject with probability tending to one, despite the null
hypothesis being correct.

To compare the variability of selected reports across studies with
that expected under prespecified Gaussian analyses, consider the
normalized between-study variation, measured by the sum of squared deviations of
the reported conclusions from their average. We normalize this quantity
by the variance \(\sigma^2/n\) of a prespecified unit-norm Gaussian
candidate and define
\[
Q_{n,K}
=
\frac{n}{\sigma^2}
\sum_{k=1}^K
\left(T_{n,k}^\star-\overline T_{n,K}^\star\right)^2.
\]

For comparison, consider a prespecified unit-norm Gaussian candidate,
represented by a fixed vector \(v\in\mathbb R^{p_n}\) with
\(\|v\|_2=1\). In study \(k\), its report is
\(T_{n,k}^{\mathrm{pre}}
=\theta+(\sigma/\sqrt n)\langle G_n^{(k)},v\rangle
\sim N(\theta,\sigma^2/n)\).
Define \(Q_{n,K}^{\mathrm{pre}}\) by replacing
\(T_{n,k}^\star\) with \(T_{n,k}^{\mathrm{pre}}\)
in the definition of \(Q_{n,K}\). Since the \(T_{n,k}^{\mathrm{pre}}\) are independent \(N(\theta,\sigma^2/n)\) variables,
\[
Q_{n,K}^{\mathrm{pre}}\sim\chi^2_{K-1},
\qquad
\E Q_{n,K}^{\mathrm{pre}}=K-1.
\]
This benchmark includes the prespecified column means of the
minimal Gaussian model in Section~\ref{sec:minimal}.
For the selected reports, however, Theorem~\ref{thm:gaussian-search}
gives
\[
\E Q_{n,K}
=
(K-1)\frac{n\operatorname{Var}(T_n^\star)}{\sigma^2}
\leqslant K-1.
\]
Search can therefore produce less between-study variation than the
prespecified Gaussian reference model predicts. In particular, a value
of \(Q_{n,K}\) in the lower tail of the reference
\(\chi^2_{K-1}\) distribution indicates unusually low between-study variability relative to that benchmark.
It does not identify selection as the cause: genuinely precise,
target-consistent estimators may also exhibit small between-study
variation. Conversely, selection need not produce under-dispersion.
Moreover, \(Q_{n,K}\) measures variability around the average reported
conclusion, not its distance from \(\theta\). For the minimal Gaussian model, Section~\ref{app:f-calibration}~of the Supplementary Material derives the corresponding exact \(F\) reference calibration for the prespecified benchmark when
\(\sigma^2\) is estimated from within-column sample variances.

The strength of this under-dispersion depends on search geometry and scaling.
Writing \(r_n=n\operatorname{Var}(T_n^\star)/\sigma^2\), one has
\(r_n\to0\) for independent Gaussian candidates as \(L_n\to\infty\),
\(r_n=1-2/\pi\) for the hypercube, \(r_n\to1/2\) for the sphere
as \(p_n\to\infty\), and \(r_n\to\rho\) under equicorrelation
when \(L_n\to\infty\) and \(\rho_n\to\rho\).
In the proxy model, an order-one displacement can even coexist with
\(r_n\to1\), for example when \(\rho_n\to0\) while
\(\rho_n\sqrt{2\log L_n/n}\to\beta>0\)
(note that this requires \(\log L_n/n\to\infty\)).
Conversely, in the independent model, \(r_n\to0\)
can coexist with \(T_n^\star\xrightarrow{P}\theta\) when
\(L_n\to\infty\) and \(\log L_n=o(n)\).
Thus under-dispersion is neither necessary nor
sufficient for a stable pseudotruth. The \(\chi^2_{K-1}\) calibration remains a benchmark for prespecified
Gaussian analyses, not a null distribution for selected reports.
In applications, between-study heterogeneity and differences in
sample sizes or search spaces may further complicate this comparison.

The Gaussian decomposition also suggests an idealized correction. Suppose
\(w_n=O(\sqrt n)\), \(\widehat\sigma_n\xrightarrow{P}\sigma\), and
\[
\frac{\widehat w_n-w_n}{\sqrt n}\xrightarrow{P}0.
\]
Then
\[
\widetilde T_n
:=
T_n^\star-\frac{\widehat\sigma_n\widehat w_n}{\sqrt n}
\xrightarrow{P}\theta.
\]
When \(w_n^2\) is of order \(n\), relative consistency
\(\widehat w_n/w_n\xrightarrow{P}1\) is sufficient. Implementing this correction requires knowledge of the search family
and the joint distribution of its candidate statistics, not merely
the selected pipeline. For a prespecified family, one possible approach
is to estimate the covariance structure and simulate the corresponding
Gaussian process, or to use a joint bootstrap when valid \citep{chernozhukov2013,white2000}. Related methods for inference on selected winners are developed by \citet{andrews2024}. This correction is model-dependent and pointwise; it makes no uniform
claim and does not contradict the impossibility results of
\citet{leebpotscher2005,leebpotscher2006,leebpotscher2008}.

Several limitations are deliberate. Richer procedures may generate pipelines
sequentially, react to intermediate results,
stop adaptively, or randomize; their final choice can be represented by
\(Q_n(\cdot\mid D_n)\), but that representation suppresses the internal search
dynamics. The correlated theory is Gaussian, while the non-Gaussian extension
treats independent candidates. The joint treatment of non-Gaussian dependence, sequential generation
of the search space, and independent replication lies beyond the scope
of the present analysis.

Replication remains essential evidence, but one must specify what is being
replicated. When the candidate analysis is fixed independently of the realized
data, agreement across independent samples has its usual sampling
interpretation. When the complete procedure searches anew in each study,
agreement instead describes the stability of the search-and-report rule, not
the validity of any particular winning analysis. Whether that stable output is close to the scientific target cannot
be established from replication agreement alone.

\bibliographystyle{apalike}
\bibliography{references_v16}

\begin{thebibliography}{}

\bibitem[Aczel et~al., 2026]{aczel2026}
Aczel, B., Szaszi, B., Clelland, H.~T., et~al. (2026).
\newblock Investigating the analytical robustness of the social and behavioural
  sciences.
\newblock {\em Nature}, 652:135--142.

\bibitem[Amelunxen et~al., 2014]{amelunxen2014}
Amelunxen, D., Lotz, M., McCoy, M.~B., and Tropp, J.~A. (2014).
\newblock Living on the edge: Phase transitions in convex programs with random
  data.
\newblock {\em Inf. Inference}, 3:224--294.

\bibitem[Andrews et~al., 2024]{andrews2024}
Andrews, I., Kitagawa, T., and McCloskey, A. (2024).
\newblock Inference on winners.
\newblock {\em Q. J. Econ.}, 139:305--358.

\bibitem[Bachoc et~al., 2020]{bachoc2020}
Bachoc, F., Preinerstorfer, D., and Steinberger, L. (2020).
\newblock Uniformly valid confidence intervals post-model-selection.
\newblock {\em Ann. Statist.}, 48:440--463.

\bibitem[Berk et~al., 2013]{berk2013}
Berk, R., Brown, L., Buja, A., Zhang, K., and Zhao, L. (2013).
\newblock Valid post-selection inference.
\newblock {\em Ann. Statist.}, 41:802--837.

\bibitem[Bertran et~al., 2026]{bertran2026}
Bertran, M., Fogliato, R., and Wu, Z.~S. (2026).
\newblock Many {AI} analysts, one dataset: Navigating the agentic data science
  multiverse.
\newblock {\em Proc. Natl. Acad. Sci. USA}, 123:e2606495123.

\bibitem[Borell, 1975]{borell1975}
Borell, C. (1975).
\newblock The {Brunn--Minkowski} inequality in {Gauss} space.
\newblock {\em Invent. Math.}, 30:207--216.

\bibitem[Botvinik-Nezer et~al., 2020]{botvinik2020}
Botvinik-Nezer, R., Holzmeister, F., Camerer, C.~F., et~al. (2020).
\newblock Variability in the analysis of a single neuroimaging dataset by many
  teams.
\newblock {\em Nature}, 582:84--88.

\bibitem[Boucheron et~al., 2013]{boucheron2013}
Boucheron, S., Lugosi, G., and Massart, P. (2013).
\newblock {\em Concentration Inequalities: A Nonasymptotic Theory of
  Independence}.
\newblock Oxford University Press, Oxford.

\bibitem[Breznau et~al., 2022]{breznau2022}
Breznau, N., Rinke, E.~M., Wuttke, A., et~al. (2022).
\newblock Observing many researchers using the same data and hypothesis reveals
  a hidden universe of uncertainty.
\newblock {\em Proc. Natl. Acad. Sci. USA}, 119:e2203150119.

\bibitem[Chatterjee, 2014]{chatterjee2014}
Chatterjee, S. (2014).
\newblock {\em Superconcentration and Related Topics}.
\newblock Springer, Cham.

\bibitem[Chernozhukov et~al., 2013]{chernozhukov2013}
Chernozhukov, V., Chetverikov, D., and Kato, K. (2013).
\newblock Gaussian approximations and multiplier bootstrap for maxima of sums
  of high-dimensional random vectors.
\newblock {\em Ann. Statist.}, 41:2786--2819.

\bibitem[Dembo and Zeitouni, 1998]{dembozeitouni1998}
Dembo, A. and Zeitouni, O. (1998).
\newblock {\em Large Deviations Techniques and Applications}.
\newblock Springer, New York, second edition.

\bibitem[Dudley, 1967]{dudley1967}
Dudley, R.~M. (1967).
\newblock The sizes of compact subsets of {Hilbert} space and continuity of
  {Gaussian} processes.
\newblock {\em J. Funct. Anal.}, 1:290--330.

\bibitem[Dwork et~al., 2015a]{dwork2015a}
Dwork, C., Feldman, V., Hardt, M., Pitassi, T., Reingold, O., and Roth, A.
  (2015a).
\newblock Preserving statistical validity in adaptive data analysis.
\newblock In {\em Proc.\ 47th ACM STOC}, pages 117--126. Association for
  Computing Machinery.

\bibitem[Dwork et~al., 2015b]{dwork2015b}
Dwork, C., Feldman, V., Hardt, M., Pitassi, T., Reingold, O., and Roth, A.
  (2015b).
\newblock The reusable holdout: Preserving validity in adaptive data analysis.
\newblock {\em Science}, 349:636--638.

\bibitem[Fithian et~al., 2014]{fithian2014}
Fithian, W., Sun, D., and Taylor, J. (2014).
\newblock Optimal inference after model selection.
\newblock {\em arXiv preprint arXiv:1410.2597}.

\bibitem[Gelman and Loken, 2013]{gelmanloken2013}
Gelman, A. and Loken, E. (2013).
\newblock The garden of forking paths: Why multiple comparisons can be a
  problem, even when there is no ``fishing expedition'' or ``p-hacking'' and
  the research hypothesis was posited ahead of time.
\newblock Unpublished manuscript.

\bibitem[Impagliazzo et~al., 2022]{impagliazzo2022}
Impagliazzo, R., Lei, R., Pitassi, T., and Sorrell, J. (2022).
\newblock Reproducibility in learning.
\newblock In {\em Proc.\ 54th ACM STOC}, pages 818--831. Association for
  Computing Machinery.

\bibitem[Kuchibhotla et~al., 2022]{kuchibhotla2022}
Kuchibhotla, A.~K., Kolassa, J.~E., and Kuffner, T.~A. (2022).
\newblock Post-selection inference.
\newblock {\em Annu. Rev. Stat. Appl.}, 9:505--527.

\bibitem[Leadbetter et~al., 1983]{leadbetter1983}
Leadbetter, M.~R., Lindgren, G., and Rootz{\'e}n, H. (1983).
\newblock {\em Extremes and Related Properties of Random Sequences and
  Processes}.
\newblock Springer, New York.

\bibitem[Ledoux and Talagrand, 1991]{ledoux1991}
Ledoux, M. and Talagrand, M. (1991).
\newblock {\em Probability in Banach Spaces: Isoperimetry and Processes}.
\newblock Springer, Berlin.

\bibitem[Leeb and P\"otscher, 2005]{leebpotscher2005}
Leeb, H. and P\"otscher, B.~M. (2005).
\newblock Model selection and inference: Facts and fiction.
\newblock {\em Econometric Theory}, 21:21--59.

\bibitem[Leeb and P\"otscher, 2006]{leebpotscher2006}
Leeb, H. and P\"otscher, B.~M. (2006).
\newblock Can one estimate the conditional distribution of post-model-selection
  estimators?
\newblock {\em Ann. Statist.}, 34:2554--2591.

\bibitem[Leeb and P\"otscher, 2008]{leebpotscher2008}
Leeb, H. and P\"otscher, B.~M. (2008).
\newblock Can one estimate the unconditional distribution of
  post-model-selection estimators?
\newblock {\em Econometric Theory}, 24:338--376.

\bibitem[Lin, 2013]{lin2013}
Lin, W. (2013).
\newblock Agnostic notes on regression adjustments to experimental data:
  Reexamining {F}reedman's critique.
\newblock {\em Ann. Appl. Stat.}, 7:295--318.

\bibitem[McLatchie and Vehtari, 2024]{mclatchie2024}
McLatchie, Y. and Vehtari, A. (2024).
\newblock Efficient estimation and correction of selection-induced bias with
  order statistics.
\newblock {\em Stat. Comput.}, 34:132.

\bibitem[{National Academies of Sciences, Engineering, and Medicine},
  2019]{nasem2019}
{National Academies of Sciences, Engineering, and Medicine} (2019).
\newblock {\em Reproducibility and Replicability in Science}.
\newblock The National Academies Press, Washington, DC.

\bibitem[Pickands, 1968]{pickands1968}
Pickands, III, J. (1968).
\newblock Moment convergence of sample extremes.
\newblock {\em Ann. Math. Statist.}, 39:881--889.

\bibitem[Resnick, 1987]{resnick1987}
Resnick, S.~I. (1987).
\newblock {\em Extreme Values, Regular Variation and Point Processes}.
\newblock Springer, New York.

\bibitem[Russo and Zou, 2016]{russo2016}
Russo, D. and Zou, J. (2016).
\newblock Controlling bias in adaptive data analysis using information theory.
\newblock In Gretton, A. and Robert, C.~C., editors, {\em Proc.\ 19th AISTATS},
  volume~51 of {\em Proceedings of Machine Learning Research}, pages
  1232--1240. PMLR.

\bibitem[Silberzahn et~al., 2018]{silberzahn2018}
Silberzahn, R., Uhlmann, E.~L., Martin, D.~P., et~al. (2018).
\newblock Many analysts, one data set: Making transparent how variations in
  analytic choices affect results.
\newblock {\em Adv. Methods Pract. Psychol. Sci.}, 1:337--356.

\bibitem[Simmons et~al., 2011]{simmons2011}
Simmons, J.~P., Nelson, L.~D., and Simonsohn, U. (2011).
\newblock False-positive psychology: Undisclosed flexibility in data collection
  and analysis allows presenting anything as significant.
\newblock {\em Psychol. Sci.}, 22:1359--1366.

\bibitem[Steegen et~al., 2016]{steegen2016}
Steegen, S., Tuerlinckx, F., Gelman, A., and Vanpaemel, W. (2016).
\newblock Increasing transparency through a multiverse analysis.
\newblock {\em Perspect. Psychol. Sci.}, 11:702--712.

\bibitem[Talagrand, 2014]{talagrand2014}
Talagrand, M. (2014).
\newblock {\em Upper and Lower Bounds for Stochastic Processes: Modern Methods
  and Classical Problems}.
\newblock Springer, Berlin.

\bibitem[Taylor and Tibshirani, 2015]{taylor2015}
Taylor, J. and Tibshirani, R.~J. (2015).
\newblock Statistical learning and selective inference.
\newblock {\em Proc. Natl. Acad. Sci. USA}, 112:7629--7634.

\bibitem[Tsirelson et~al., 1976]{tsirelson1976}
Tsirelson, B.~S., Ibragimov, I.~A., and Sudakov, V.~N. (1976).
\newblock Norms of {Gaussian} sample functions.
\newblock In Maruyama, G. and Prokhorov, J.~V., editors, {\em Proc.\ Third
  Japan--USSR Symp.\ Probab.\ Theory}, volume 550 of {\em Lecture Notes in
  Math.}, pages 20--41. Springer, Berlin.

\bibitem[Vershynin, 2018]{vershynin2018}
Vershynin, R. (2018).
\newblock {\em High-Dimensional Probability: An Introduction with Applications
  in Data Science}.
\newblock Cambridge University Press, Cambridge.

\bibitem[White, 2000]{white2000}
White, H. (2000).
\newblock A reality check for data snooping.
\newblock {\em Econometrica}, 68:1097--1126.

\bibitem[Yu and Kumbier, 2020]{yukumbier2020}
Yu, B. and Kumbier, K. (2020).
\newblock Veridical data science.
\newblock {\em Proc. Natl. Acad. Sci. USA}, 117:3920--3929.

\end{thebibliography}

\clearpage
\appendix
\section{Supplementary Material}

\subsection{Proof of Proposition~\ref{prop:exponential-search}}
\label{app:proof-exponential}

For every \(L\geqslant1\), the map
\(z\mapsto\max_{1\leqslant j\leqslant L}z_j\) is \(1\)-Lipschitz.
The Gaussian Poincar\'e inequality \citep{boucheron2013} therefore gives
\(\operatorname{Var}(M_L)\leqslant1\). For two independent repetitions,
\[
\E\!\left[\bigl(T_n^{\star(1)}-T_n^{\star(2)}\bigr)^2\right]
=\frac{2\sigma^2}{n}\operatorname{Var}(M_{L_n})
\leqslant\frac{2\sigma^2}{n},
\]
which proves \eqref{eq:10} without any restriction on \(L_n\).
Symmetry and almost-sure uniqueness make \(J_n^\star\) uniform on
\(\{1,\ldots,L_n\}\); independence of the two experiments gives
\eqref{eq:11}.

Under the additional assumption
\(\log L_n/n\to c\), if \(c>0\), then \(L_n\to\infty\), and \eqref{eq:8} gives
\[
\frac{M_{L_n}}{\sqrt n}
=
\frac{M_{L_n}}{\sqrt{2\log L_n}}
\sqrt{\frac{2\log L_n}{n}}
\xrightarrow{P}\sqrt{2c},
\]
which proves \eqref{eq:9}. If \(c=0\), for every \(\varepsilon>0\),
\[
\mathbb P(M_{L_n}>\varepsilon\sqrt n)
\leqslant
L_n\mathbb P(Z>\varepsilon\sqrt n)
\leqslant
\exp\!\left(\log L_n-\frac{\varepsilon^2n}{2}\right),
\qquad Z\sim N(0,1),
\]
which tends to zero, while
\(\mathbb P(M_{L_n}<-\varepsilon\sqrt n)\leqslant
\mathbb P(Z_{n,1}<-\varepsilon\sqrt n)\to0\). Thus \eqref{eq:9} also holds for
\(c=0\). 

\subsection{Variance of the Gaussian maximum}
\label{app:gaussian-maximum-variance}
Let
\[
a_L=\sqrt{2\log L},
\qquad
b_L=a_L-\frac{\log\log L+\log(4\pi)}{2a_L}.
\]
Classical extreme-value theory gives, as \(L\to\infty\),
\(a_L(M_L-b_L)\Rightarrow \Lambda\), where \(\Lambda\) is standard Gumbel with
\(\operatorname{Var}(\Lambda)=\pi^2/6\), and convergence of the corresponding second
moments \citep{pickands1968} yields
\[
\operatorname{Var}(M_L)
\sim
\frac{\pi^2}{6a_L^2}
=
\frac{\pi^2}{12\log L}.
\]
See \citet{leadbetter1983,resnick1987}.

\subsection{Proof of the threshold-stopping claims}
\label{app:threshold}

Let \(q_n=\mathbb P(T_{n,1}\geqslant t)
=\overline\Phi(\sqrt n(t-\theta)/\sigma)\), where \(t>\theta\)
and \(\overline\Phi=1-\Phi\). Gaussian tail asymptotics give
\[
\frac{1}{n}\log q_n\longrightarrow
-\frac{(t-\theta)^2}{2\sigma^2}.
\]
Under the stated assumption \(\log L_n/n\to c>(t-\theta)^2/(2\sigma^2)\),
we have \(L_nq_n\to\infty\). Consequently, the probability that no
candidate crosses the threshold is given by
\[
(1-q_n)^{L_n}\leqslant \exp(-L_nq_n)\longrightarrow0.
\]
Given at least one crossing, the first reported value has the
law of \(T_{n,1}\) conditional on \(T_{n,1}\geqslant t\), by independence
of the candidates. For every \(\varepsilon>0\),
\[
\mathbb P(T_{n,1}\geqslant t+\varepsilon\mid T_{n,1}\geqslant t)
=
\frac{\overline\Phi(\sqrt n(t+\varepsilon-\theta)/\sigma)}
{\overline\Phi(\sqrt n(t-\theta)/\sigma)}
\longrightarrow0.
\]
Since the no-crossing event is asymptotically negligible, the reported estimator
converges in probability to \(t\).

Let \(J_n=0\) if no candidate crosses the threshold
(the fallback outcome), and let \(J_n=j\) if candidate \(j\) is the first to cross it.
Then
\[
\mathbb P(J_n=0)=(1-q_n)^{L_n},
\]
whereas, for \(1\leqslant j\leqslant L_n\),
\[
\mathbb P(J_n=j)=(1-q_n)^{j-1}q_n.
\]
Consequently,
\[
\max_{0\leqslant j\leqslant L_n}\mathbb P(J_n=j)
\leqslant \max\{q_n,(1-q_n)^{L_n}\}
\longrightarrow 0.
\]
Thus, for two independent repetitions,
\[
\begin{aligned}
\mathbb P(J_n^{(1)}=J_n^{(2)})
&=\sum_{j=0}^{L_n}\mathbb P(J_n=j)^2\\
&\leqslant \max_{0\leqslant j\leqslant L_n}
\mathbb P(J_n=j)
\longrightarrow 0.
\end{aligned}
\]

\subsection{Proof of Proposition~\ref{prop:proxy-search}}
\label{app:proof-proxy}

The maximizing index depends only on
\((S_{n,1},\ldots,S_{n,L_n})\). Since the \(\eta_{n,j}\)'s are independent of
all selection scores, \(\eta_{n,J_n^\star}\) is standard Gaussian and
independent of \(M_{L_n}\), which gives \eqref{eq:proxy-decomposition};
\eqref{eq:proxy-variance} follows immediately. Moreover,
\[
\frac{\rho_nM_{L_n}}{\sqrt n}
=
\rho_n\sqrt{\frac{2\log L_n}{n}}
\frac{M_{L_n}}{\sqrt{2\log L_n}}
\xrightarrow{P}\beta,
\]
while \(\sqrt{1-\rho_n^2}E_n/\sqrt n\xrightarrow{P}0\). This proves
\eqref{eq:proxy-limit}.

\subsection{Proof of Proposition~\ref{prop:hypercube}}
\label{app:hypercube-proof}

From \eqref{eq:38},
\[
w(V_p)
=
\frac1{\sqrt p}\sum_{k=1}^p\E|G_k|
=
\sqrt{\frac{2p}{\pi}},
\]
which gives \eqref{eq:39}. Independence and
\(\operatorname{Var}(|G_1|)=1-2/\pi\) give \eqref{eq:40}. For \(p=p_n\),
\[
T_n^\star
=
\theta+\frac{\sigma}{\sqrt{np_n}}
\sum_{k=1}^{p_n}|G_{n,k}|.
\]
Hence
\[
\E T_n^\star-\theta
=
\sigma\sqrt{\frac{2p_n}{\pi n}},
\qquad
\operatorname{Var}(T_n^\star)
=
\frac{\sigma^2}{n}\left(1-\frac{2}{\pi}\right).
\]
For two independent repetitions,
the exact variance and independence of the
copies give \eqref{eq:42}.
The maximizing signs are
\(\operatorname{sgn}(G_{n,1}),\ldots,
\operatorname{sgn}(G_{n,p_n})\),
which are independent Rademacher variables.
Hence the maximizing vertex is uniform on
the \(2^{p_n}\) vertices, giving
\eqref{eq:43}.
If, in addition, \(p_n/n\to\kappa\in[0,\infty)\),
the mean converges to
\(\theta+\sigma\sqrt{2\kappa/\pi}\),
while the variance tends to zero,
proving \eqref{eq:41}.

\subsection{Calculations for the spherical example}
\label{app:sphere-calculations}

For \(G\sim N(0,I_p)\), \(\|G\|_2\) has a chi distribution with \(p\) degrees
of freedom and
\[
\E\|G\|_2
=
\sqrt2\,\frac{\Gamma((p+1)/2)}{\Gamma(p/2)}
=
\sqrt p\left(1-\frac{1}{4p}+O(p^{-2})\right).
\]
Hence
\[
w(S^{p-1})^2=p-\frac12+O(p^{-1}),
\qquad
\operatorname{Var}(\|G\|_2)=p-w(S^{p-1})^2\longrightarrow\frac12.
\]
If \(A=G/\|G\|_2\), symmetry gives \(\E A=0\) and
\(\E A_jA_k=0\) for \(j\neq k\), while exchangeability and
\(\|A\|_2=1\) give \(\E A_j^2=1/p\). For independent copies
\(A^{(1)},A^{(2)}\), this yields
\[
\E\langle A^{(1)},A^{(2)}\rangle=0,
\qquad
\E\langle A^{(1)},A^{(2)}\rangle^2=\frac1p.
\]

\subsection{Proof of Corollary~\ref{cor:geometric-regimes}}
\label{app:entropy-corollary}

If \(\mathcal J_n=o(\sqrt n)\), Proposition~\ref{prop:entropy-bounds} gives
\(w_n=o(\sqrt n)\). Borell--TIS concentration then gives
\(T_n^\star=\theta+\sigma w_n/\sqrt n+O_P(n^{-1/2})\), and hence target
consistency. Under \eqref{eq:60}, Sudakov's lower bound at
\(\varepsilon_0\) yields
\[
w_n\geqslant c_0\varepsilon_0
\sqrt{\log M(\mathcal A_n,d_n,\varepsilon_0)},
\]
which gives \eqref{eq:61} and then \eqref{eq:62}. Finally, applying the same Gaussian concentration bound to each of two
independent repetitions gives
\(T_n^{\star(r)}=\theta+\sigma w_n/\sqrt n+O_P(n^{-1/2})\),
\(r=1,2\), and therefore
\[
T_n^{\star(1)}-T_n^{\star(2)}\xrightarrow{P}0.
\]

\subsection{Proof of Proposition~\ref{prop:large-deviation-search}}
\label{app:proof-large-deviation}

It suffices to consider \(\varepsilon>0\) such that
\(t_c-\varepsilon>\theta\) and \(t_c+\varepsilon<t_+\). By the union bound and
\eqref{eq:70}--\eqref{eq:72},
\[
\mathbb P(T_n^\star\geqslant t_c+\varepsilon)
\leqslant
L_n\mathbb P(T_{n,1}\geqslant t_c+\varepsilon)\longrightarrow0,
\]
because the logarithmic rate is
\(c-I(t_c+\varepsilon)<0\). For the lower tail, with
\(q_n(\varepsilon)=\mathbb P(T_{n,1}\geqslant t_c-\varepsilon)\),
\[
\mathbb P(T_n^\star<t_c-\varepsilon)
=
(1-q_n(\varepsilon))^{L_n}
\leqslant
\exp\{-L_nq_n(\varepsilon)\},
\]
and
\(n^{-1}\log(L_nq_n(\varepsilon))\to c-I(t_c-\varepsilon)>0\). Hence the lower-tail probability also vanishes, proving \eqref{eq:73}.
The conclusion-replicability statement follows immediately. If the candidate
distribution is continuous, symmetry makes the unique maximizing index uniform.
Therefore, for two independent repetitions,
\[
\mathbb P(J_n^{\star(1)}=J_n^{\star(2)})
=
\sum_{j=1}^{L_n}\frac1{L_n^2}
=
\frac1{L_n}.
\]

\subsection{Exact pooled-variance calibration in the minimal Gaussian model}
\label{app:f-calibration}

Assume \(n\geqslant2\). Recall the minimal Gaussian model of Section~\ref{sec:minimal}. In each of \(K\) independent
studies, we observe \(L_n\) columns of \(n\) observations each, with
all \(Y_{ij}^{(k)}\) independent and distributed as
\(N(\theta,\sigma^2)\). Study \(k\) selects the column
\(J_{n,k}^\star\) with the largest empirical mean and reports
\(T_{n,k}^\star=\overline Y_{\cdot,J_{n,k}^\star}^{(k)}\).
Let \(S_{n,k}^{\star2}\) denote the usual sample variance of the
observations in that selected column.

The vector of column sample variances is independent of the vector of column means. Since the
column variances are identically distributed and selection depends only on the
means, the selected sample variance has the same conditional law given all
column means as any fixed column variance, and is independent of those means.
The same argument applies jointly across the independent studies. Hence
\[
\frac{(n-1)S_{n,k}^{\star2}}{\sigma^2}\sim\chi^2_{n-1}
\]
even after selection. With
\(\widehat\sigma_n^2=K^{-1}\sum_{k=1}^KS_{n,k}^{\star2}\),
\[
\frac{K(n-1)\widehat\sigma_n^2}{\sigma^2}\sim\chi^2_{K(n-1)}.
\]
Replacing the unknown variance \(\sigma^2\) in \(Q_{n,K}\)
by its pooled estimator gives
\[
\widehat Q_{n,K}
=
\frac{\sigma^2}{\widehat\sigma_n^2}Q_{n,K}
=
\frac{n}{\widehat\sigma_n^2}
\sum_{k=1}^K
(T_{n,k}^\star-\overline T_{n,K}^\star)^2.
\]
For comparison, let \(T_{n,k}^{\mathrm{pre}}\) denote the mean of a
column fixed in advance in study \(k\), and define
\[
\widehat Q_{n,K}^{\mathrm{pre}}
=
\frac{n}{\widehat\sigma_n^2}
\sum_{k=1}^K
\left(T_{n,k}^{\mathrm{pre}}
-\overline T_{n,K}^{\mathrm{pre}}\right)^2.
\]
As established above, the pooled estimator \(\widehat\sigma_n^2\) is independent of all
column means, hence of \(Q_{n,K}^{\mathrm{pre}}\). Consequently,
\[
\frac{\widehat Q_{n,K}^{\mathrm{pre}}}{K-1}
\sim F_{K-1,K(n-1)}.
\]
This is an exact reference calibration, not the distribution of
\(\widehat Q_{n,K}/(K-1)\) after selection.

\end{document}